\documentclass[aap]{imsart}

\RequirePackage{amsthm,amsmath,amsfonts,amssymb}
\RequirePackage[numbers]{natbib}
\usepackage[utf8x,utf8]{inputenc}
\usepackage{amsfonts, amsmath, amsthm, array, csquotes, dsfont}
\usepackage{color}
\usepackage{algorithm}
\usepackage[noend]{algpseudocode}
\usepackage{url}

\startlocaldefs
\theoremstyle{plain}
\newtheorem{lemma}{Lemma}[section]
\newtheorem{corollary}[lemma]{Corollary}
\newtheorem{thm}[lemma]{Theorem}
\newtheorem{prop}[lemma]{Proposition}
\theoremstyle{remark}
\newtheorem{rk}[lemma]{Remark}

\newtheorem*{HI-det}{Deterministic induction hypothesis}
\newtheorem*{HI-stoch}{Stochastic induction hypothesis}
\newtheorem{defn}[lemma]{Definition}
\renewcommand{\d}{\textrm{d}}

\newcommand{\E}{\mathbb{E}}
\newcommand{\Indicator}[1]{\mathds{1}\left(#1\right)}

\renewcommand{\P}{\mathbb{P}}
\newcommand{\R}{\mathbb{R}}
\newcommand{\tr}{\text{\textnormal{tr}}}
\newcommand{\Var}{\mathbb{V}\text{\textnormal{ar}}}

\newcommand{\lone}[1]{\lvert #1 \rvert}
\newcommand{\Lone}[1]{\left \lvert #1 \right \rvert}

\newcommand{\proj}{{\textnormal{\textrm{proj}}}}
\newcommand{\tar}{{\textnormal{\textrm{tar}}}}

\newcommand{\frob}[1]{\lVert #1 \rVert}

\newcommand{\N}{{\mathbb N}}

\newcommand{\calE}{{\cal E}}

\newcommand{\calM}{{\cal M}}

\newcommand{\calS}{{\cal S}}

\newcommand{\grad}{\nabla}

\newcommand{\nfin}{{n_p}}
\newcommand{\npar}{{n_g}}
\endlocaldefs

\begin{document}

\begin{frontmatter}
\title{Insight from the Kullback--Leibler divergence into adaptive importance sampling schemes for rare event analysis in high dimension}
%\title{A sample article title with some additional note\thanksref{T1}}
\runtitle{Adaptive importance sampling in high dimension}
%\thankstext{T1}{A sample of additional note to the title.}

\begin{aug}
%%%%%%%%%%%%%%%%%%%%%%%%%%%%%%%%%%%%%%%%%%%%%%%
%% Only one address is permitted per author. %%
%% Only division, organization and e-mail is %%
%% included in the address.                  %%
%% Additional information can be included in %%
\author[A]{\fnms{Jason}~\snm{Beh}\ead[label=e1]{jason.beh@onera.fr}},
\author[B]{\fnms{Yonatan}~\snm{Shadmi}\ead[label=e2]{shadmi@campus.technion.ac.il}}
\and
\author[C]{\fnms{Florian}~\snm{Simatos}\ead[label=e3]{florian.simatos@isae-supaero.fr}} %C doesnt work
%%%%%%%%%%%%%%%%%%%%%%%%%%%%%%%%%%%%%%%%%%%%%%
%% Addresses                                %%
%%%%%%%%%%%%%%%%%%%%%%%%%%%%%%%%%%%%%%%%%%%%%%
\address[A]{ONERA/DTIS, Université de Toulouse, F-31055 Toulouse, France\printead[presep={,\ }]{e1}}
\address[B]{Technion\printead[presep={,\ }]{e2}}
\address[C]{Fédération ENAC ISAE-SUPAERO ONERA, Université de Toulouse, France\printead[presep={,\ }]{e3}}
\end{aug}

\begin{abstract}
	We study two adaptive importance sampling schemes for estimating the probability of a rare event in the high-dimensional regime $d \to \infty$ with $d$ the dimension. The first scheme is the prominent cross-entropy (CE) method, and the second scheme, motivated by recent results, uses as auxiliary distribution a projection of the optimal auxiliary distribution on a lower dimensional subspace. In these schemes, two samples are used: the first one to learn the auxiliary distribution and the second one, drawn according to the learned distribution, to perform the final probability estimation. Contrary to the common belief that the sample size needs to grow exponentially in the dimension to make the estimator consistent and avoid the weight degeneracy phenomenon, we find that a polynomial sample size in the first learning step is enough. We prove this result assuming that the sought probability is bounded away from~$0$. For CE, insight is provided on the polynomial growth rate which remains implicit. In contrast, we study the second scheme in a simple computational framework assuming that samples from the conditional distribution are available. This makes it possible to show that the sample size only needs to grow like $rd$ with $r$ the effective dimension of the projection, which highlights the potential benefits of these projection methods.
\end{abstract}

\begin{keyword}[class=MSC]
\kwd[Primary ]{65C05}
\kwd{62F12}
\kwd[; secondary ]{62-08}
\end{keyword}

\begin{keyword}
\kwd{Importance sampling; curse of dimensionality; weight degeneracy; Cross entropy; dimension reduction;}
\end{keyword}

\end{frontmatter}
%%%%%%%%%%%%%%%%%%%%%%%%%%%%%%%%%%%%%%%%%%%%%%
%% Please use \tableofcontents for articles %%
%% with 50 pages and more                   %%
%%%%%%%%%%%%%%%%%%%%%%%%%%%%%%%%%%%%%%%%%%%%%%
\tableofcontents

%%%%%%%%%%%%%%%%%%%%%%%%%%%%%%%%%%%%%%%%%%%%%%
%%%% Main text entry area:
\section{Introduction}

We are interested in the problem of estimating by importance sampling (IS) the probability $p_f(A) = \P_f(X \in A)$ in a high-dimensional setting: $f$ is some distribution on $\R^d$, $X$ a random variable taking values in $\R^d$ with distribution $f$, and $A$ a measurable subset of $\R^d$, with $d$ large. Given an auxiliary importance sampling density $g$, the IS estimator of $p_f(A)$ is given by
\begin{equation} \label{eq:p}
	\hat p_f(A) = \frac{1}{\nfin} \sum_{i=1}^{\nfin} \frac{f(Y_i)}{g(Y_i)} \xi_A(Y_i) \text{ with the $Y_i$'s i.i.d. $\sim g$}
\end{equation}
where here and in the sequel, $\xi_B$ is the characteristic function of the set $B$: $\xi_B(x) = \Indicator{x \in B}$. The curse-of-dimensionality stipulates that, as the dimension $d$ of the problem grows large, the sample size $\nfin$ needs to grow exponentially in $d$, see for instance~\cite{Au03:0, Bengtsson08:0}. Underlying this statement is the fact that we are interested in the high-dimensional regime $d \to \infty$ and all other quantities implicitly depend on $d$: however, this dependency is omitted from the notation for the sake of simplicity. For instance, the sample size $\nfin$ is chosen as a function of $d$ ($\nfin = \nfin(d)$), and the curse-of-dimensionality asserts that one needs to take $\nfin \gg e^{\alpha d}$ for some $\alpha > 0$ in order to make $\hat p_f(A)$ consistent, meaning for instance that $\hat p_f(A) / p_f(A) \to 1$ in some suitable sense. Here and in the sequel, for two sequences $u,v$ implicitly or explicitly depending on $d$, the notation $u \gg v$ means that $u/v \to \infty$ as $d \to \infty$.

In this paper, we study IS schemes where the auxiliary density is estimated: we have a target auxiliary density $g_\tar$ which is estimated by $\hat g_\tar$ using a sample of size $\npar$. These schemes belong to the generic family of adaptive importance sampling methods, which involve three steps: sampling, weighting and an adaptation step, which can be iterated. Most of the popular IS schemes, such as population monte-carlo or adaptive multiple importance sampling, are actually adaptive importance sampling schemes~\cite{Bugallo17:0}. In these schemes, it is important to keep in mind that there are two distinct sample sizes: $\npar$ is the sample size used to estimate the target auxiliary density $g_\tar$ in the adaptation step, whereas $\nfin$ refers to the IS sample size used to estimate $p_f(A)$ with $g_\tar$ as in~\eqref{eq:p}. Our main result is that the curse-of-dimensionality can be avoided provided $\npar$ only grows polynomially in $d$. More precisely, we study two different families of target auxiliary densities $g_\tar$ and for each, we show that if $\npar$ grows polynomially in $d$ (with the exponent depending on the given target auxiliary density), then $\nfin$ does not need to grow exponentially in $d$: actually, any sequence $\nfin \to \infty$ growing to infinity makes $\hat p_f(A)$ consistent. Said otherwise, our results show that the curse-of-dimensionality can be avoided in adaptive importance sampling schemes, provided the sample size of the adaptation step grows polynomially, and not exponentially, in the dimension.

Our results also shed light on the weight-degeneracy phenomenon, which states that, as the dimension increases, the largest importance sampling weight takes all the mass. One way to formulate the weight-degeneracy phenomenon is that, as $d \to \infty$, we have
\[ \frac{\max_{i=1, \ldots, \nfin} (f/g)(Y_i)}{\sum_{i=1, \ldots, \nfin} (f/g)(Y_i)} \Rightarrow 1 \]
with $\Rightarrow$ denoting convergence in distribution. Such a behavior clearly prevents importance sampling estimators to converge, and this is why a large literature has been devoted to avoiding this phenomenon (see the literature overview in Section~\ref{sub:lit-overview}). Moreover, Chatterjee and Diaconis have recently proposed to use this ratio for testing for convergence~\cite[Section 2]{Chatterjee18:0}. Our results show at the same time that the importance sampling estimator $\hat p_f(A)$ is consistent, and that weight degeneracy is avoided.

The paper is organized as follows. Section~\ref{sec:main-results} presents our main results, Theorems~\ref{thm:CE},~\ref{thm:g_proj} and~\ref{thm:p}, together with a literature overview and an overview of the proofs. Section~\ref{sec:preliminary-results} gathers preliminary results used in the proofs of the main theorems. We then begin with the easier proofs of Theorems~\ref{thm:CE} and~\ref{thm:p} in Sections~\ref{sec:proof-g_proj} and~\ref{sec:proof-p}, respectively, and then provide the proof of Theorem~\ref{thm:CE} in Section~\ref{sec:proof-CE}. Finally, the proof of a technical proposition is provided in the appendix.

\section{Main results} \label{sec:main-results}

We will study three different target distributions $g_\tar$, namely:
\begin{itemize}
    \item $g_A$ is the optimal Gaussian density in the sense of the Kullback--Leibler divergence;
    \item $g_\proj$ is obtained by projecting the covariance matrix of $g_A$ onto a lower-dimensional subspace;
    \item $g_t$ is the auxiliary density in the $t$-th step of the cross-entropy (CE) method.
\end{itemize}
Note that CE is an algorithm tailored for estimating small probabilities, while the projection methods underlying $g_\proj$ have been proposed to make reliability methods work in high dimension. For each of these auxiliary densities, we are interested in whether $\hat p_f(A)$ is consistent, and also whether weight-degeneracy is avoided. To capture both properties, we define the notion of high-dimensional efficiency. In the following definition, the distribution $g$ may be random: then $\E(\phi(Y) \mid g)$ with $Y \sim g$ is a notation for
\[ \E \left( \phi(Y) \mid g \right) = \int_{\R^d} \phi(y) g(y) \d y. \]
Actually, $g$ will be a Gaussian density with some random parameter $(\mu, \Sigma)$, and so conditioning on $g$ is tantamount to conditioning on $(\mu, \Sigma)$.

\begin{defn} [High-dimensional efficiency for $A$]
For each dimension $d$, let $A$, $f$, $g$ and the $Y_i$'s be as in~\eqref{eq:p} (with $g$ potentially random), and let in addition $\ell = f/g$.

As $d \to \infty$, we say that the sequence of auxiliary distributions $g$ is efficient in high dimension for $A$ if, for any sequence $\nfin \to \infty$, the two following conditions hold:
	\begin{equation} \label{eq:HDE-WD}
		\E \left ( \frac{\max_{i=1, \ldots, \nfin} \ell(Y_i)}{\sum_{i=1, \ldots, \nfin} \ell(Y_i)} \mid g \right) \Rightarrow 0
	\end{equation}
	and
	\begin{equation} \label{eq:HDE-L1}
		\E \left ( \left \lvert \frac{1}{p_f(A) \nfin} \sum_{i=1}^\nfin \ell(Y_i) \xi_A(Y_i) - 1 \right \rvert \mid g \right) \Rightarrow 0.
	\end{equation}
\end{defn}

What is important in this definition is that the sampling size $\nfin$ does not need to grow at some prescribed rate with the dimension: thus, this avoids the curse-of-dimensionality in a strong sense. Chatterjee and Diaconis~\cite{Chatterjee18:0} proved that the minimal sampling size for an IS scheme is of the order of $e^{D(f||g)}$ or $e^{D(f|_A||g)}$ with $f|_A$ the distribution $f$ conditioned on $A$ and $D(h||g)$ the Kullback--Leibler divergence between two densities $h$ and $g$. That the sampling size may grow at any speed actually hinges upon the fact that $D(f||g)$ and $D(f|_A||g)$ remain bounded, which is the kind of results that we will prove in this paper.

Note moreover that only~\eqref{eq:HDE-L1} depends on $A$, but the idea is that $g$ will be chosen as a function of $A$, which makes~\eqref{eq:HDE-WD} implicitly depend on $A$ as well. As will be seen below, the price to pay will be in the sampling size in the adaptation step where the auxiliary density is learned (in particular, $g$ will be taken as an estimator $\hat g_\tar$ of some target density $g_\tar$), but in this step, the sampling size will only need to grow polynomially in the dimension, and not exponentially as when the curse-of-dimensionality occurs.

Our results rely on two main assumptions: that $X$ follows a standard Gaussian distribution, and that $p_f(A)$ as a sequence in $d$ is bounded away from $0$. The Gaussian assumption is mostly made for technical reasons but covers a large range of applications since general transformations allow to express a random variable $X$ as $\Phi(Y)$ with $Y$ Gaussian~\cite{lebrun09b,lebrun09c}. When such transformations can be applied, then we can re-express the probability as $\P(X \in A) = \P(Y \in \Phi^{-1}(A))$ and then apply our results with $A' = \Phi^{-1}(A)$.  Moreover, the Gaussian assumption is sometimes encountered in practical settings, see for instance the examples from~\cite{Ehre23:0} and~\cite{Uribe20:0} discussed in Section~\ref{sub:discussion-p}. On the other hand, the assumption $\inf_d p_f(A) > 0$ is discussed in detail in Section~\ref{sub:discussion-p}.

\subsection{Minimal notation} \label{sub:min-notation}

We first introduce the minimal set of notation necessary in order to state our main results, while further notation is introduced in Section~\ref{sub:further-notation}. Let in the sequel $\calS_d$ denote the space of $d \times d$ symmetric, semi-definite positive matrices. For $\mu \in \R^d$ and $\Sigma \in \calS_d$, we denote by $N(\mu, \Sigma)$ the $d$-dimensional Gaussian distribution with mean $\mu$ and covariance matrix $\Sigma$. In the rest of the paper, we consider the case where the initial distribution $f$ is the density of a $d$-dimensional standard Gaussian vector in dimension~$d$, i.e.,
\[ f(x) = (2\pi)^{-d/2} e^{-\lVert x \rVert^2 / 2}, \ x \in \R^d, \]
where here and in the sequel, $\lVert x \rVert$ denotes the $L_2$-norm of some vector $x \in \R^d$ (note also that here and elsewhere, we identify a distribution with its density). For any density $g$ on $\R^d$ and any measurable set $B \subset \R^d$, we denote by $p_g(B) = \int \xi_B g$ the measure of the set $B$ under~$g$, and $g|_B = g \xi_B / p_g(B)$ the distribution $g$ conditioned on $B$. Concerning random variables, we will adopt the following convention:
\begin{itemize}
 \item $X$ will refer to a generic random variable, and its distribution will be indicated by a subscript in the probability or expectation: for instance, $\E_f(X)$ is the mean of $X$ under $\P_f$, i.e., when $X$'s distribution is $f$;
 \item we will use $Y$ to refer to random variables drawn according to a given distribution: in this case, their mean will be denoted by the generic $\E$.
\end{itemize}
For instance, when the $Y_i$'s are i.i.d.\ drawn according to $g$, then we will write
\[ \E \left( \frac{1}{n} \sum_{i=1}^n \frac{f(Y_i)}{g(Y_i)} \phi(Y_i) \right) = \E\left(\frac{f(Y_i)}{g(Y_i)}\phi(Y_i)\right) = \E_g\left(\frac{f(X)}{g(X)}\phi(X)\right) = \E_f(\phi(X)). \]
Another example is the probability $p_g(B)$ which can equivalently be written as $p_g(B) = \P_g(X \in B) = \P(Y \in B)$ with $Y \sim g$.

For $x \in \R_+$, we denote by $[x] = \max\{n \in \N: n \leq x\}$ its integer part. For a c\`adl\`ag function $f: \R_+ \to \R_+$, we consider its left-continuous inverse
\[ f^{-1}(t) = \inf\{s \geq 0: f(s) \geq t\}, \ t \geq 0. \]
Note that for $x, t \geq 0$ we have $f(x) \geq t \Leftrightarrow x \geq f^{-1}(t)$ and that if $f$ is continuous, then $f(f^{-1}(t)) = t$.

\subsection{High-dimensional efficiency of CE densities} \label{sub:CE}

We begin with results on the CE densities $\hat g_t$. Our approach is to first study a deterministic version of CE corresponding to target densities $g_t$, and then to study the true version of CE, which corresponds to using as auxiliary density an estimation $\hat g_t$ of $g_t$.

The CE method works for sets $A$ of the form $A = \{x \in \R^d: \varphi(x) \geq q\}$ for some measurable function $\varphi:\R^d \to \R$ and some threshold $q \in \R$. Algorithm~\ref{alg:CE-det} presents the deterministic version of the CE method, leading to the target densities $g_t$. This is the algorithm one would implement if all quantities were known and need not be estimated. However, in practice it is intractable because it relies on the quantiles $q_t$ (Step 2a) and on the conditional mean $\mu_t$ and variance $\Sigma_t$ under $f|_{A_t}$ (Step 2b) which cannot be computed analytically. This leads to the true version of CE described in Algorithm~\ref{alg:CE}, where the untractable quantities $q_t$, $\mu_t$ and $\Sigma_t$ are replaced by estimation. Note that Algorithm~\ref{alg:CE} uses another sequence $m$ which is the size of the sample used in the quantile estimation step 2a. As for the other sequences, $m = m(d)$ is implicitly a sequence depending on $d$. Note also that there is usually a stopping criterion in step $2$, typically when $q_t \geq q$, i.e., $A_t \subset A$. Here we do not consider the stopping criterion, and identify conditions under which $g_t$ and $\hat g_t$ are efficient in high dimension for $A$. In the following statement and in the sequel, we say that a function $\varphi: \R^d \to \R$ has no atom if for every $x \in \R$, the set $\varphi^{-1}(\{x\})$ has zero Lebesgue measure.

\begin{algorithm}[t]
\caption{Deterministic version of CE}\label{alg:CE-det}
\begin{algorithmic}
\Require $\rho \in (0,1)$
\State 1. define $\mu_0 = 0$ and $\Sigma_0 = I$ and start with $g_0 = f = N(\mu_0, \Sigma_0)$;
\State 2. Iterate the following steps:
\State \ \ \ \ (a) given $g_t$, consider $q_t = F^{-1}(1-\rho)$, where $F$ is the cumulative distribution function of $\varphi(X)$ under $\P_{g_t}$;
\State \ \ \ \ (b) define $A_t = \{x: \varphi(x) > q_t\}$ and
 	\[ \mu_{t+1} = \mu_{A_t} = \E_{f|_{A_t}}(X) \ \text{ and } \ \Sigma_{t+1} = \Sigma_{A_t} = \Var_{f|_{A_t}}(X); \]
  \State \ \ \ \ (c) define $g_{t+1} = N(\mu_{t+1}, \Sigma_{t+1})$, let $t = t+1$ and go to Step~(a).
\end{algorithmic}
\end{algorithm}

\begin{algorithm}[t]
\caption{Stochastic version of CE}\label{alg:CE}
\begin{algorithmic}
\Require $\rho \in (0,1)$, sample sizes $\npar$ and $m$
\State 1. start with $\hat g_0 = f$;
\State 2. Iterate the following steps:
\State \ \ \ \ (a) given $\hat g_t$, draw $Y'_1, \ldots, Y'_m$ i.i.d.\ according to $\hat g_t$, rank them in increasing values of $\varphi$: $\varphi(Y'_{(1)}) \leq \cdots \leq \varphi(Y'_{(m)})$ and let $\hat q_t = \varphi(Y'_{([(1-\rho) m])})$;
\State \ \ \ \ (b) define $\hat A_t = \{x: \varphi(x) > \hat q_t\}$, $\ell = f/\hat g_t$, draw $Y_1, \ldots, Y_\npar$ i.i.d.\ according to $\hat g_t$, independently from the $Y'_k$'s, and compute 
 	\begin{equation} \label{eq:CE-p-mu}
			\hat p_t = \frac{1}{\npar} \sum_{i=1}^\npar \ell(Y_i) \xi_{\hat A_t}(Y_i), \ \hat \mu_{t+1} = \frac{1}{\npar \hat p_t} \sum_{i=1}^\npar \ell(Y_i) \xi_{\hat A_t}(Y_i) Y_i
		\end{equation}
 	and
 	\begin{equation} \label{eq:CE-Sigma}
			\hat \Sigma_{t+1} = \frac{1}{\npar \hat p_t} \sum_{i=1}^\npar \ell(Y_i) \xi_{\hat A_t}(Y_i) Y_i Y_i^\top - \hat \mu_{t+1} \hat \mu_{t+1}^\top;
		\end{equation}
  \State \ \ \ \ (c) define $\hat g_{t+1} = N(\hat \mu_{t+1}, \hat \Sigma_{t+1})$, let $t = t+1$ and go to Step~(a).
\end{algorithmic}
\end{algorithm}

\begin{thm}\label{thm:CE}
    Assume that:
    \begin{itemize}
        \item $\inf_d p_f(A) > 0$;
        \item $\inf_d \rho > 0$;
        \item for every $d$, $\varphi$ has no atom.
    \end{itemize}
    Then for every $t \geq 0$, $g_t$ is efficient in high dimension for $A$.\\
    Assume in addition that $m \to \infty$. Then for each $t \geq 0$, there exists a finite constant $\kappa_t > 0$ such that if $\npar \gg d^{\kappa_t}$, then $\hat g_t$ is efficient in high dimension.
\end{thm}

\begin{rk} \label{rk:dep}
    In true CE schemes, the same sample is used in steps 2a and 2b to estimate both the quantile and the mean and variance. The simpler dependency structure in Algorithm~\ref{alg:CE} makes the theoretical analysis easier. We leave it as open research to study the algorithms with the full dependency structure.
\end{rk}

\begin{rk}
    We study here the adaptive-level version of CE, but the tools developed herein could also be used to study the fixed-level version of CE, close to the subset simulation algorithm~\cite{Au01:0, Cerou12:0, Douc07:0}.
\end{rk}

\noindent \textit{Discussion on the constants $\kappa_t$.} Let us now discuss the constant $\kappa_t$. For $t = 0$ we have $\kappa_0 = 1$, and for $t \geq 1$, we are only able to prove the existence of some $\kappa_t > 0$. To give some intuition on this constant, let us introduce the notation $\lambda_1(\Sigma)$ for the smallest eigenvalue of a symmetric, positive definitive matrix $\Sigma$. Let further $\hat \lambda_{*,t} = \min\{\lambda_1(\hat \Sigma_1), \ldots, \lambda_1(\hat \Sigma_t)\}$ and $\hat \kappa_{*,t} = 8 \max(1, 1/\hat \lambda_{*,t}-1)$, so that $\hat \kappa_{*,t} = 8$ if $\lambda_1(\hat \Sigma_k) \geq 1/2$ for $k = 1, \ldots, t$, and $\hat \kappa_{*,t} = 8 (1/\hat \lambda_{*,t} - 1) > 8$ otherwise. In Section~\ref{subsub:proof-outline} below, we explain that if $\npar \gg d^\kappa$ for some $\kappa > \hat \kappa_{*,t}$, then we could prove that $\hat g_t$ is efficient in high dimension for $A$. This would give a more explicit expression for the exponent of the required growth rate, but this would not be satisfactory because the growth rate would be random.

As $\hat \Sigma_t$ is an estimator of $\Sigma_t$, it is clear that this result suggests that Theorem~\ref{thm:CE} should hold for any $\kappa_t > \kappa_{*,t}$ with $\kappa_{*,t} = 8 \max(1, 1/\lambda_{*,t}-1)$ with $\lambda_{*,t} = \min\{\lambda_1(\Sigma_1), \ldots, \lambda_1(\Sigma_t)\}$. Because of monotonicity, in order to establish such a result, it would be enough to prove that
\begin{equation} \label{eq:kappa}
	\forall \varepsilon > 0,\ \P(\lambda_1(\hat \Sigma_t) \geq (1-\varepsilon) \lambda_1(\Sigma_t)) \to 1.
\end{equation}
However, it is well-known that controlling the smallest eigenvalue of random matrices is a difficult task, see for instance~\cite{Bai93-0}, and we did not manage to find simple arguments to prove~\eqref{eq:kappa}. However, we managed to prove the existence of some $\underline \lambda_t > 0$ such that $\P(\lambda_1(\hat \Sigma_t) \geq \underline \lambda_t) \to 1$, and then Theorem~\ref{thm:CE} holds with $\kappa_t = 8 \max(1, 1/\underline \lambda_{*,t}-1)$ with $\underline \lambda_{*,t} = \min\{\underline \lambda_1, \ldots, \underline \lambda_t\}$. We believe that, upon additional technical assumptions (e.g., on the growth rate of $m$ and regularly properties for $\varphi$), one could prove something like~\eqref{eq:kappa} and therefore relate $\kappa_t$ to the $\lambda_1(\Sigma_t)$'s. However, our main objective was to show that polynomial growth rates were enough, and so we content ourselves with the result as stated above, although it could most probably be strengthened along various directions.

Note finally that the reason why smallest eigenvalues play a role in our proofs is that we need finite $\alpha$-th moments of the likelihood ratios $f/\hat g_t$. More precisely, we need $\alpha > 0$ such that
\[ \E_f \left[ \left( \frac{f(X)}{\hat g_t(X)} \right)^\alpha \mid \hat g_t \right] < \infty, \ \text{ almost surely}. \]
But for this to hold, one needs
\[ \alpha < \min \left( 1, \frac{\lambda_1(\hat \Sigma_t)}{1-\lambda_1(\hat \Sigma_t)} \right) \]
which is where the smallest eigenvalues kick in.

\subsection{High-dimensional efficiency of auxiliary distributions using projection methods in a simple computational framework} \label{sub:proj}

It is well-known that the optimal IS auxiliary density for estimating the probability $p_f(A)$ is $f|_A = f \xi_A / p_f(A)$, i.e., the distribution $f$ conditioned on~$A$. Indeed, for this choice of auxiliary density, we have $\hat p = p_f(A)$ (with $\hat p$ defined in~\eqref{eq:p}), i.e., $p_f(A)$ is perfectly estimated. Of course, $f|_A$ is intractable as it involves the unknown quantity $p_f(A)$. Among Gaussian auxiliary densities, the one that minimizes the Kullback--Leibler divergence with $f|_A$ is $g_A = N(\mu_A, \Sigma_A)$ with $\mu_A$ and $\Sigma_A$ the mean and variance of~$f|_A$:
\[ \mu_A = \E_{f|_A}(X) \ \text{ and } \ \Sigma_A = \Var_{f|_A}(X) = \E_{f|_A}(XX^\top) - \mu_A \mu_A^\top \]
which makes $g_A$ a natural candidate for a good auxiliary density. Actually, most if not all AIS schemes are somehow geared toward this density~\cite{Bugallo17:0}. In CE for instance, the $g_t$'s are meant to get close to $g_A$, so that $\hat g_t$ may be thought of as an approximation of $g_A$.

When the dimension is large, several authors have proposed targeting other densities $g_\proj$ obtained by projecting $\Sigma_A$ onto a low-dimensional subspace. Various subspaces on which to project were proposed recently~\cite{El-Masri21:0, El-Masri:0, Uribe20:0}, and they all lead to considering a Gaussian auxiliary density $g_\proj = N(\mu_A, \Sigma_\proj)$ with mean $\mu_A$ and variance $\Sigma_\proj$ defined as
\begin{equation} \label{eq:Sigma-proj}
	\Sigma_\proj = \sum_{k=1}^r (v_k - 1) d_k d_k^\top + I \ \text{ with } \ v_k = d_k^\top \Sigma_A d_k
\end{equation}
where the $d_i$'s form an orthonormal family, and $r$ is the dimension of the small subspace on which to project. In practice, we have $r \leq 3$ most of the times, but our results will apply to any $r \leq d$. They apply in particular for $r = d$ in which case we have $g_\proj = g_A$, and so $g_A$ can be seen as special case of $g_\proj$. Several choices are considered in~\cite{El-Masri21:0, El-Masri:0, Uribe20:0}:
\begin{itemize}
	\item in~\cite{Uribe20:0}, a smooth approximation $\tilde \xi_A \approx \xi_A$ of the characteristic function is considered. The $d_k$'s are the eigenvectors of the matrix $H := \E_{f|_A}( (\grad \log \tilde \xi_A(X)) (\grad \log \tilde \xi_A(X))^\top )$ and they are ranked in decreasing order of the corresponding eigenvalues, i.e., $d_1$ corresponds to the largest eigenvalue of $H$, $d_2$ to the second largest eigenvalue, etc;
	\item in~\cite{El-Masri21:0}, only one direction is considered ($r = 1$) and $d_1 = \mu_A / \lVert \mu_A \rVert$;
	\item in~\cite{El-Masri:0}, the $d_k$'s are the eigenvectors of $\Sigma_A$, and they are ranked in decreasing order according to the image by the function $h(x) = x-1- \log x$ of the eigenvalues, i.e., $d_1$ is associated to the eigenvalue maximizing $h$, etc.
\end{itemize}
These different choices were analyzed in~\cite{El-Masri21:0, El-Masri:0, Uribe20:0} and were found to perform very well numerically. However, an analytic explanation for that success was not complete and this work makes a step in this direction.

Similarly as for CE, the target densities $g_\proj$ are intractable because they involve the unknown conditional mean and variance $\mu_A$ and $\Sigma_A$. Here we place ourselves in a simple computational framework, and we assume that samples from $f|_A$ are available from which $\mu_A$ and $\Sigma_A$ can be estimated, see the discussion following Corollary~\ref{cor:g_A}. We therefore assume that we are given a sample $(Y_{A,i})_i$ of i.i.d.\ random variables drawn according to $f|_A$, and this sample is used to estimate $g_A$ and $g_\proj$ as follows: $\hat g_A = N(\hat \mu_A, \hat \Sigma_A)$ with
\begin{equation} \label{eq:hat-mu-A-hat-Sigma-A}
	\hat \mu_A = \frac{1}{\npar} \sum_{k=1}^\npar Y_{A,k} \ \text{ and } \ \hat \Sigma_A = \frac{1}{\npar} \sum_{k=1}^\npar Y_{A,k} Y_{A,k}^\top - \hat \mu_A \hat \mu_A^\top
\end{equation}
and $\hat g_\proj = N(\hat \mu_A, \hat \Sigma_\proj)$ with
\begin{equation} \label{eq:Sigma-proj-hat}
	\hat \Sigma_\proj = \sum_{k=1}^r (\hat v_k - 1) d_k d_k^\top + I \ \text{ with } \ \hat v_k = d_k^\top \hat \Sigma_A d_k
\end{equation}
where the $d_k$'s form an orthonormal family. The $d_k$'s are allowed to be random, and they should typically be thought of as estimators of deterministic target projection directions. In our results, the $d_k$'s are assumed to be independent from the $Y_{A,i}$'s, see Remark~\ref{rk:dep} below.

The following result is in the same vein as Theorem~\ref{thm:CE} for CE, although simpler: the target density $g_\proj$ is efficient in high dimension for $A$ under the condition $\inf_d p_f(A) > 0$, and an additional growth rate condition on $\npar$ is needed for the same to hold for $\hat g_\proj$.

\begin{thm}\label{thm:g_proj}
	Assume that $\inf_d p_f(A) > 0$, and consider any $r \leq d$ and any orthonormal family $(d_1, \ldots, d_r)$. Then $g_\proj$ is efficient in high dimension for $A$.\\
    If in addition $\npar \gg rd$ and $(d_1, \ldots, d_r)$ is independent from the $Y_{A,i}$'s, then $\hat g_\proj$ is efficient in high dimension for $A$.
\end{thm}

It is striking to note that the only influence of the projection family $(d_1, \ldots, d_r)$ is through its size $r$. In practice however, the choice of the projection family has a significant impact on the quality of the final estimation. Further results are thus called upon in order to capture the influence of the $d_k$'s.

Since $g_\proj = g_A$ when $r = d$, we readily obtain the following result.

\begin{corollary} \label{cor:g_A}
	Assume that $\inf_d p_f(A) > 0$. Then $g_A$ is efficient in high dimension for $A$.\\
 If moreover $\npar \gg d^2$, then $\hat g_A$ is efficient in high dimension for $A$.
\end{corollary}

\begin{rk}
 The algorithms studied here slightly differ from those proposed in the literature because of different dependency structures. More precisely, in~\cite{El-Masri21:0, El-Masri:0, Uribe20:0}, the directions $d_k$ on which to project are indeed estimations of deterministic target directions, but in practice they are computed from $\hat \Sigma_A$ and are thus not independent from the $Y_{A,i}$'s.
\end{rk}

\noindent \textit{Discussion on the simple computational framework.} In Theorem~\ref{thm:g_proj} we assume that samples from $f|_A$ are available, which is admittedly a strong assumption that limits the applicability of this result. Still, this simple computational framework is insightful for various reasons.

First, since $f|_A$ is known up to a normalizing constant, various simulation schemes (typically, MCMC) can be used to sample from it and so in simple settings, this assumption can be satisfied.

Second, our results give valuable insight into the high-dimensional behavior of IS. Indeed, $g_A$ is the target distribution of most AIS schemes and so results on $\hat g_A$ in a favorable computational context provide potentially useful benchmark results. For instance, already comparing Corollary~\ref{cor:g_A} and Theorem~\ref{thm:g_proj} provides useful insight into the potential benefits of projection schemes; likewise, comparing Theorem~\ref{thm:CE} with Corollary~\ref{cor:g_A} highlights the influence of having to use IS methods in the estimations of $\mu_t$ and $\Sigma_t$. We leave it as open question to study mixed schemes $g_{t, \proj}$ relying on both CE and projection as for instance in~\cite{Uribe20:0}.
\\

\noindent \textit{Discussion on $\hat g_A$ vs $\hat g_\proj$.} Let us discuss some insight in $\hat g_A$ and $\hat g_\proj$ provided by Theorem~\ref{thm:g_proj}. If one could sample from $g_A$ and $g_\proj$, even though both are efficient in high dimension according to Theorem~\ref{thm:g_proj} and its Corollary~\ref{cor:g_A}, it is clear that $g_A$ would be preferable since it is the optimal auxiliary density. However, $g_\proj$ involves fewer parameters and is therefore intuitively easier to estimate. Thus, although $g_A$ is better than $g_\proj$, $\hat g_A$ incurs more estimation error which could make $\hat g_\proj$ preferable to $\hat g_A$. Theorem~\ref{thm:g_proj} and Corollary~\ref{cor:g_A} provide evidence in that direction, in that they show that $\hat g_A$ remains efficient in high dimension provided $\npar \gg d^2$, whereas for $\hat g_\proj$, one only needs $\npar \gg r d$. As mentioned earlier, in practice we typically have $r \leq 3$, and so one only needs a linear growth rate for $\hat g_\proj$, but a quadratic growth rate for $\hat g_A$.

Of course these results do not claim that these growth rates are sharp, and that the conditions $\npar \gg d^2$ and $\npar \gg rd$ are necessary for $\hat g_A$ and $\hat g_\proj$ to be efficient in high dimension. Nonetheless, the following result suggests that the $d^2$ threshold is sharp. In the following result, we assume that $\mu_A$ and $\Sigma_A$ are estimated from a sample $(Y_{A,k})$ drawn according to $N(\mu_A, \Sigma_A)$ instead of $f|_A$: since by definition, $N(\mu_A, \Sigma_A)$ and $f|_A$ have the same mean and variance, drawing the $Y_{A,k}$'s according to $N(\mu_A, \Sigma_A)$ in~\eqref{eq:hat-mu-A-hat-Sigma-A} still gives consistent estimators. Of course this scheme is of no practical interest, as there does not seem to be methods to sample from $N(\mu_A, \Sigma_A)$ without knowing $\mu_A$ and $\Sigma_A$. However, this scheme presents a theoretical interest, in that if the $Y_{A,k}$'s are Gaussian, then $\hat \mu_A$ is Gaussian and $\hat \Sigma_A$ follows a Wishart distribution. In this case, explicit formulas are available which allows to prove the following result.
\begin{prop}\label{pro:E-DKL}
	Assume that in~\eqref{eq:hat-mu-A-hat-Sigma-A} the $Y_{A,k}$'s are i.i.d.~drawn according to $N(\mu_A, \Sigma_A)$ instead of $f|_A$. Assume that $\npar \gg d$: then $\sup_d \E(D(f || \hat g_A)) < \infty$ if $\npar \gg d^2$, and $\E(D(f || \hat g_A)) \to \infty$ if $\npar \ll d^2$.
\end{prop}

The proof of this result is given in the appendix. As mentioned previously in the introduction, Chatterjee and Diaconis~\cite{Chatterjee18:0} proved that the sampling size needs to be at least $e^{D(f||\hat g_A)}$ in order for the IS estimator to be close to its target. Thus, the fact that the expected KL divergence diverges for $\npar \ll d^2$ is an indication that the $d^2$ threshold is sharp, in that if $\npar \ll d^2$, then there is a minimal growth rate imposed upon $\nfin$, namely $e^{D(f||\hat g_A)}$, and so $\hat g_A$ cannot be efficient in high dimension, at least in the way we defined it.\\

\subsection{Discussion of the assumption $\inf_d p_f(A) > 0$} \label{sub:discussion-p}

Here we discuss the assumption $\inf_d p_f(A) > 0$ required for all our results. First, note that although the regime $p_f(A) \to 0$ seems quite natural, we are not aware of previous work studying this regime in the same context as in the present paper, even in fixed dimension. In particular, the assumption $\inf_d p_f(A) > 0$ is enforced in previous works that study high-dimensional importance sampling in a reliability context, see for instance~\cite{Au03:0}. This assumption is relevant or even required in different cases. For instance, Bassamboo et al.~\cite{Bassamboo08} use IS to estimate the probability of large portfolio losses in some high-dimensional regime: there the dimension represents the number of ``obligors'', and the authors consider both cases where the sought probability goes to $0$ or stays away from $0$ as the dimension increases.

There are also more general contexts in which the assumption $\inf_d p_f(A) > 0$ is relevant. One example is when $p_f(A)$ represents the probability stemming from an approximation scheme. In this case, it is natural to expect $p_f(A) \approx p_{\text{true}}$ with $p_{\text{true}}$ the probability of interest for the ``true'' system: provided $p_{\text{true}} > 0$, it is then natural to expect that $\inf_d p_f(A) > 0$. To illustrate this idea, let us consider two engineering examples from~\cite[Section~$4.3.1$]{Ehre23:0} and~\cite[Section~$5.3$]{Uribe20:0}. There the authors are interested in some probability $p_{\text{true}}$ of the form $p_{\text{true}} = \P(u^*(\mathcal{E}) \in U)$ with $u^*(\mathcal{E})$ the solution to a partial differential equation (PDE) whose coefficients are driven by a Gaussian field $(\mathcal{E}(y), y \in D)$ with $D$ the domain of the PDE. To give a concrete example, Uribe et al.\ consider the following equation:
\[ G(y) \nabla^2 u(y) + \frac{\mathcal{E}(y)}{2(1-\nu)} \nabla (\nabla \cdot u(y)) + b = 0, \]
see~\cite[Section~$5.3$]{Uribe20:0} for more details. Then $u^*(\mathcal{E})$ would denote the solution to this PDE, with an emphasis on the dependency of the solution to the random coefficient $\mathcal{E}$. The Karhunen-Loève expansion amounts to approximating the Gaussian field $\mathcal{E}$ by
\[ \mathcal{E}(y) \approx E_d(y, X) = \exp \left( \mu + \sigma \sum_{i=1}^d \sqrt{\lambda_i} \phi_i(y) X_i \right) \]
where $X = (X_1, \ldots, X_d)$ is a standard Gaussian vector. The other coefficients $\mu, \sigma, \gamma_i$ and $\phi_i$ are deterministic but beyond this, they do not play a role for the present discussion, so the interested reader is for instance referred to~\cite{Ehre23:0},~\cite{Uribe20:0} or~\cite{betz14} for more details. The key point is that we have an approximation of the form $\mathcal{E} \approx E_d(X)$ with $d$ the number of terms in the Karhunen-Loève expansion and $X$ a standard Gaussian vector which drives the randomness since $E_d$ is a deterministic function. Plugging in this approximation into $p_{\text{true}}$, this suggests the approximation
\[ p_{\text{true}} = \P(u^*(\mathcal{E}) \in U) \approx \P(u^*(E_d(X)) \in U). \]
Defining $A = \{x: u^*(E_d(x)) \in U\}$, this leads to the approximation $p_{\text{true}} \approx p_f(A)$, under which the relevant regime is indeed the regime where $\inf_d p_f(A) > 0$\footnote{More precisely, the approximation $p_{\text{true}} \approx p_f(A)$ would translate in mathematical terms as $p_f(A) \to p_{\text{true}}$, which implies that $\liminf p_f(A) > 0$ provided $p_{\text{true}} > 0$. At this point note that: 1/ since all our results are in the asymptotic regime $d \to \infty$, they continue to hold if $\liminf p_f(A) > 0$ instead of $\inf p_f(A) > 0$ and 2/ the assumptions $\liminf_d p_f(A) > 0$ and $\inf p_f(A) > 0$ only differ when $p_f(A) = 0$ for some $d$, which can be ruled out in most practical applications, for instance the example discussed here.}. Note that in this example, the dimension $d$ represents the number of terms in the Karhunen-Loève expansion.

One insight of our work is that if $\inf_d p_f(A) > 0$, then $D(f|_A||g_A)$ and $D(f||g_A)$ are bounded (see Corollary~\ref{cor:Sigma-proj-mu-A-bounded} and the proof of Proposition~\ref{prop:first-step}). As mentioned earlier, the results of Chatterjee and Diaconis~\cite{Chatterjee18:0} suggest that there is no minimal growth rate for $\nfin$. The following result shows that if $p_f(A) \to 0$, then either $D(f|_A||g_A)$ or $D(f||g_A)$ is unbounded, which imposes a minimal growth rate on $\nfin$. In the naive Monte-Carlo scheme, the required number of samples grow like $1/p_f(A)$, while in Section~\ref{sec:proof-p} we will prove that $D(f|_A||g_A) \leq -\log p_f(A)$, suggesting that IS may require less sample than MC. Further investigation on minimal growth rates for $\nfin$ when $p_f(A) \to 0$ represents an interesting research question which we leave untouched, and here we content ourselves with the following result.

\begin{thm}\label{thm:p}
	Assume that the condition $p_f(A) \to 0$ holds. Then we have either $\sup_d D(f || g_A) = \infty$ or $\sup_d D(f|_A || g_A) = \infty$.
\end{thm}

\subsection{Literature overview} \label{sub:lit-overview}

\subsubsection{Importance sampling as a sampling scheme}

Importance sampling is a popular numerical method that can be used for sampling from an intricate distribution and for reducing the variance of a Monte-Carlo estimator, see for instance~\cite{owen_monte_2013, Robert04:0} for a general introduction. The literature on the former case of using IS for sampling is very large. Starting from the basic IS schemes, many improved variants have been proposed, using mixtures and control variates~\cite{Owen00-0}, resampling schemes~\cite{Rubin87-0, Rubin87-1}, use of particular auxiliary densities~\cite{Hesterberg95-0} or local MCMC-like moves~\cite{neal2001}, to list only a few. Moreover, instead of aiming to sample from a given distribution, one may instead aim to sample from a sequence of distributions, leading to so-called sequential MC or IS schemes, see for instance~\cite{Del-Moral06:0, Doucet01-1}. Sequential schemes can also be used in static contexts~\cite{Chopin02-0}, and this idea lead to the fundamental population Monte-Carlo algorithm and its variants~\cite{Cappe08:0, Cappe04:0}. Finally, adaptive IS schemes involve learning steps whereby parameters of the auxiliary distribution are updated against past performance~\cite{Bugallo17:0, Oh92:0}.

The theoretical analysis of the basic IS scheme is straightforward: as a sum of i.i.d.\ random variables, its consistency is settled by the law of large numbers and its speed by the central limit theorem. However, in more advanced schemes, resampling and recycling of samples create intricate dependency structures which make the analysis challenging. Theoretical results on the behavior of complex adaptive IS schemes can for instance be found in~\cite{Akyildiz21:0, Douc07:0, Douc07-0, marin_consistency_2019, Portier18:0}.

Concerning the high-dimensional behavior of IS, \enquote{it is well known that importance sampling is usually inefficient in high-dimensional spaces}~\cite{Doucet01-1}. One of the main reasons is the weight degeneracy problem, whereby the largest IS weight takes all the mass. This phenomenon is due to the potential heavy tail of likelihood ratios, which arises in high dimension as the densities naturally become singular with respect to one another. For this reason, various schemes have been proposed by transforming the weights in order to reduce their variance~\cite{El-Laham18:0, Elvira17:0, ionides_truncated_2008, Koblents15:0, Martino18:0, rubinstein09, vehtari2022pareto}.

Although verified empirically, to our knowledge weight degeneracy has only been addressed theoretically in the related context of particle filters~\cite{Bengtsson08:0} (see also~\cite{Surace19-0} for a review), where it is proved that the sample size needs to grow exponentially in the dimension in order to avoid weight degeneracy. In an unpublished report~\cite{Li05-0}, the same authors have additional results on IS where for an i.i.d.\ target and a test function that only depends on one coordinate, they claim that the sample size needs to grow at least in $\exp(d^{1/3})$ with $d \to \infty$ the dimension (see in particular~\cite[Proposition $3.6$]{Li05-0}\footnote{Note that in this result, the authors assume that $\E_{g_1}((f_1(X(1))/g_1(X(1)))^a) < \infty$ for some $a > 0$, where $f_1$ is the first marginal of $f$ and $g_1$ of $g$, and assuming $f$ and $g$ i.i.d.. However, as $\E_{g_1}(f_1(X(1))/g_1(X(1)) = 1$, this assumption seems to always hold with $a = 1$. We think that the finite-moment assumption probably needs to hold for some $a > 1$ rather than $a > 0$, see the related Remark~\ref{rk:WD} below.}). High-dimensional results are also obtained in~\cite{Beskos14:0, Beskos14-0}, who consider an i.i.d.\ target that is approximated by a sequential MC scheme through bridging distributions. Among other results, they prove that provided the number of bridging densities grows linearly in the dimension, the effective sample size remains bounded, therefore suggesting that \enquote{AIS schemes may beat the curse of dimensionality in some scenarios if properly designed}~\cite{Bugallo17:0}. Our main results point toward a similar conclusion in the context of rare event probability estimation.

\subsubsection{Importance sampling in a reliability context}

In reliability, the overarching theme is the estimation of the probability $p := \P_f(X \in A)$ of an important event $A$ (e.g., the failure of a critical component) which is deemed rare, so that $p$ is small. The coefficient of variation $\sqrt{\Var(\hat I)} / \E(\hat I)$ of the naive MC estimator $\hat I$ scales like $1/\sqrt{p}$, which calls upon improved techniques, such as the prominent subset simulation method~\cite{Au01:0, Cerou12:0, Cheng22-0, Rashki21-0}. 

In the IS realm, various schemes have also been proposed, see for instance~\cite{au1999, Cheng23-0, owen2019, Papaioannou16-0} and~\cite{morio2015, morio14, Tabandeh22-0} for a review. Recall from its description in Algorithms~$1$ and $2$ above that CE aims at building a sequence of densities $g_t$ getting closer and closer to the optimal density $g_A$: CE can thus be seen as a special case of sequential IS, but because of its importance for our work, we will reserve a special discussion to CE below.

In high dimension, auxiliary distributions specific to reliability problems have been proposed to avoid weight degeneracy~\cite{Chiron23-0, Papaioannou19:0, Wang16:0}. From a theoretical perspective, the authors study in~\cite{Au03:0} the performance of IS for high-dimensional reliability problems. The set-up is quite similar to ours, as the authors assume that the probability $p$ is bounded away from $0$, and they also consider Gaussian auxiliary distributions (they also consider mixtures of Gaussian distributions, but do not have theoretical results in this case). Their main result is that in order for IS to be applicable in high dimension in the case where the initial distribution is standard normal, the covariance matrix of the auxiliary density must be a finite-rank perturbation of the identity. This is very close in spirit to what we prove here, as our proofs essentially rely on proving that $\frob{\Sigma - I}$ remains bounded, with $\Sigma$ the covariance matrix of the auxiliary density considered. Note however that a significant difference between~\cite{Au03:0} and our results is that the authors in~\cite{Au03:0} consider the variance as the performance metric, which imposes a restriction on the set of auxiliary distributions that can be considered. More precisely, in order for $f(X) / g(X)$ to have a finite second moment, with $X \sim g$, $f = N(0,I)$ and $g = N(\mu, \Sigma)$, all eigenvalues of $\Sigma$ must be larger than $1/2$. If this condition is violated, then the authors conclude in~\cite{Au03:0} that IS is not applicable; since we consider the $L_1$ norm, our scheme still works. Note however that, as explained in the discussion following Theorem~\ref{thm:CE}, the threshold $1/2$ has a strong impact on the performance of CE because of its influence on the growth rate $\kappa_t$.

To conclude this literature overview, let us focus more precisely on the CE method~\cite{Boer05:0, Kroese13:0, Rubinstein04:0}. In low dimension, numerous numerical results tend to suggest that CE and its improvements are quite efficient, see for instance~\cite{chan_improved_2012, Geyer19-0, Papaioannou19:0}. However, even in this case, theoretical results backing up these numerical observations are pretty scarce. We are only aware of~\cite{Homem-de-Mello07:0} (which provides proofs of some results announced earlier in~\cite{homem-de-mello_rare_2002}) which provides theoretical guarantees on the convergence of a modified version of CE. In high dimension, CE may suffer from weight degeneracy similarly as for general IS schemes discussed above. The dimension-reduction strategies discussed above aim at avoiding this problem~\cite{El-Masri21:0, El-Masri:0, Uribe20:0}. Thus, to the best of our knowledge, our results are the first ones to theoretically address the behavior of CE in high dimension.

\subsection{Proof overview} The rest of the paper is devoted to proving Theorems~\ref{thm:CE},~\ref{thm:g_proj} and~\ref{thm:p}. Before delving into the proofs, let us give an overview. In Section~\ref{sec:preliminary-results} we gather preliminary results. This culminates in Proposition~\ref{prop:first-step}, which establishes a sufficient condition on $\mu$ and $\Sigma$ for a sequence of Gaussian densities $g = N(\mu, \Sigma)$ to be efficient in high-dimension for~$A$. Informally, the condition states that $1/p_f(A)$ and $\mu$ are bounded, that $\Sigma$ stays close to the identity, and that the spectrum of $\Sigma$ remains in a compact set of $(0,\infty)$. In Section~\ref{sec:proof-g_proj} we show that these properties hold for $g_\proj$ and $\hat g_\proj$ under appropriate conditions, and in Section~\ref{sec:proof-CE} we show that they hold for $g_t$ and $\hat g_t$.

In both cases, we first establish the efficiency in high dimension for the deterministic target densities: for $g_\proj$ this is done in Corollary~\ref{cor:Sigma-proj-mu-A-bounded}, and for $g_t$ in Proposition~\ref{prop:HDE-g_t}. Once $g_\proj$ and $g_t$ are shown to be efficient in high dimension, we then transfer the results on $\hat g_\proj$ and $\hat g_t$ by controlling the error. The most difficult part lies in controlling the covariance matrices $\hat \Sigma_\proj$ and $\hat \Sigma_t$: as explained above, one has to show that they remain close to the identity, and that their spectrum remains in a compact set of $(0,\infty)$. As we have first shown such properties to hold for $\Sigma_\proj$ and $\Sigma_t$, the idea is to transfer these properties to $\hat \Sigma_\proj$ and $\hat \Sigma_t$ by using results from random matrix theory to control the errors $\hat \Sigma_\proj - \Sigma_\proj$ and $\hat \Sigma_t - \Sigma_t$.

Finally, we note that the proofs for $g_t$ and especially $\hat g_t$ are significantly more difficult than for $g_\proj$ and $\hat g_\proj$. In particular, the proof operates by induction on $t$, and one of the key idea is to identify a suitable induction hypothesis that goes through. The deterministic induction hypothesis for $g_t$ is stated in Section~\ref{sub:g_t}, and its stochastic version for $\hat g_t$ in Section~\ref{subsub:proof-outline}.

\section{Preliminary results} \label{sec:preliminary-results}

\subsection{Further notation} \label{sub:further-notation}

Recall the notation already introduced in Section~\ref{sub:min-notation}: here, we complement this notation with further notation needed in the paper. In the sequel, a vector $x \in \R^d$ will be considered as a column vector, and its coordinates will be written $x(1), \ldots, x(d)$. For vectors, indices will refer to sequences, for instance, $X_1, \dots, X_n$ will typically denote an i.i.d.\ sequence of~$\R^d$-valued random variables drawn according to $f$, and $X_i(k)$ will denote $X_i$'s $k$-th coordinate. Let in the sequel $\calM_d$ denote the space of $d \times d$ matrices, and recall that $\calS_d$ denotes the space of $d \times d$ symmetric, semi-definite positive matrices. For a matrix $M \in \calM_d$, we will write its entries either by $M(i,j)$ or by $M_{ij}$.

For $x \in \R^d$ and $M \in \calM_d$, we denote by $\lvert x \rvert$ and $\lvert M \rvert$ the sum of the absolute values of its coordinates or entries:
\[ \lvert x \rvert = \sum_{k=1}^d \lvert x(k) \rvert \ \text{ and } \ \lvert M \rvert = \sum_{i,j} \lvert M_{ij} \rvert. \]
Note that we omit the dependency in the dimension in that $\lvert \cdot \rvert$ denotes the $L$\textsubscript{1} norm for different dimensions. This abuse of notation will be enforced throughout the paper as most of the times, dependency on $d$ will be omitted in order to ease the notation. Let further $\lVert x \rVert^2$ and $\frob{M}^2$ denote the sum of the square of its coordinates or entries:
\[ \lVert x \rVert^2 = \sum_{k=1}^d x(k)^2 \ \text{ and } \ \frob{M}^2 = \sum_{i,j} M_{ij}^2. \]
Note that $\lVert x \rVert \leq \lvert x \rvert$ and $\lVert M \rVert \leq \lvert M \rvert$. Further, for $M \in \calM_d$ a square matrix, $\frob{M}$ is its Frobenius norm, and we have $\frob{M}^2 = \tr(M M^\top)$. We denote by $\det(M)$ its determinant, and if $M$ is symmetric, we denote by $\lambda_1(M) \leq \cdots \leq \lambda_d(M)$ its eigenvalues ranked in increasing order. We will use repeatedly and without notice the variational characterization of eigenvalues, which implies in particular that, for $M \in \calM_d$ symmetric,
\[ \lambda_1(M) \lVert x \rVert^2 \leq x^\top M x \leq \lambda_d(M) \lVert x \rVert^2, \ x \in \R^d. \]

 Concerning the $L_1$ matrix norm, we will use the following result.

\begin{lemma}\label{lemma:L1-matrix}
	For $\Sigma, \Sigma' \in \calM_d$ symmetric, we have $\lvert \lambda_1(\Sigma) - \lambda_1(\Sigma') \rvert \leq \lVert \Sigma - \Sigma' \rVert$.
\end{lemma}

\begin{proof}
	Let $v \in \R^d$ with $\lVert v \rVert = 1$. Then we have
	\[ v^\top \Sigma v = v^\top (\Sigma - \Sigma') v + v^\top \Sigma' v \geq v^\top (\Sigma - \Sigma') v + \lambda_1(\Sigma'). \]
	Moreover,
	\[ \left( v^\top (\Sigma - \Sigma') v \right)^2 \leq \max \left( \lambda_1(\Sigma - \Sigma')^2, \lambda_d(\Sigma - \Sigma')^2 \right) \leq \lVert \Sigma - \Sigma' \rVert^2 \]
	and so
	\[ v^\top \Sigma v \geq \lambda_1(\Sigma') - \lVert \Sigma - \Sigma' \rVert. \]
	Since this holds for any $v \in \R^d$ with unit norm, this entails
	\[ \lambda_1(\Sigma) \geq \lambda_1(\Sigma') - \lVert \Sigma - \Sigma' \rVert \]
	which gives the result by symmetry between $\Sigma$ and $\Sigma'$.
\end{proof}

We define the function $\Psi: \calS_d \to \R_+$ by
\begin{equation} \label{eq:Psi}
	\Psi(\Sigma) = \frac{1}{2} (\tr(\Sigma) - \log \det(\Sigma) - d).
\end{equation}
Note that if $\psi(x) = x - \log x - 1$ for $x > 0$, then we have $\Psi(\Sigma) = \frac{1}{2} \sum_{i=1}^d \psi(\lambda_i(\Sigma))$. As $\psi \geq 0$ and $\psi(x) = 0 \Leftrightarrow x = 1$, this shows that $\Psi(\Sigma) \geq 0$ and that $\Psi(\Sigma) = 0 \Leftrightarrow \Sigma = I$, with $I$ the identity matrix.

Given two density functions $g$ and $g'$ on $\R^d$, with $g$ absolutely continuous with respect to $g'$, we define the Kullback--Leibler (KL) divergence between $g$ and $g'$ by
\[ D(g || g') = \int g(x) \log (g(x)/g'(x))dx = \E_g \left( \log \left( \frac{g(X)}{g'(X)} \right) \right). \]

For $B \subset \R^d$ measurable and $g$ a density on $\R^d$, we denote by $\mu^g_B$ and $\Sigma^g_B$ the mean and variance of $g|_B$:
\begin{equation} \label{eq:def-mu-sigma-g-B}
	\mu^g_B = \E_{g|_B}(X) \ \text{ and } \ \Sigma^g_B = \Var_{g|_B}(X).
\end{equation}
When $g = f$ (the standard Gaussian density), we omit the superscript and simply write
\begin{equation} \label{eq:def-mu-sigma-f-B}
	\mu_B = \mu^f_B \ \text{ and } \ \Sigma_B = \Sigma^f_B.
\end{equation}

Finally, we use $\Rightarrow$ to denote convergence in distribution, and we say that a sequence $X$ of real-valued random variables, implicitly indexed by the dimension~$d$, is bounded with high probability (whp) if there exists $K \geq 0$ such that $\P(\lvert X \rvert \leq K) \to 1$ as $d \to \infty$.

\subsection{Results from Chatterjee and Diaconis~\cite{Chatterjee18:0}} \label{sub:CD}

In order to study the high-dimensional efficiency of some sequence of auxiliary distributions, we will crucially rely on the recent results of Chatterjee and Diaconis~\cite{Chatterjee18:0}: this result shows that it is enough to focus on the KL divergence and on the tail behavior of the log-likelihood. According to Theorem~1.1 in~\cite{Chatterjee18:0}, for any measurable $\phi: \R^d \to \R$ and any $n \geq e^{D(f || g)}$, we have
\begin{multline} \label{eq:bound-CD}
 \E \left( \left \lvert \frac{1}{n} \sum_{i=1}^n \ell(Y_i) \phi(Y_i) - \E_f(\phi(X)) \right \rvert \right) \leq \left( \E_f \left(\phi(X)^2 \right) \right)^{1/2} \times\\
 \left[ \left(\frac{e^{D(f || g)}}{n}\right)^{1/4} + 2 \left( \P_f \left(L(X) \geq \frac{1}{2} \log n + \frac{1}{2}D(f||g) \right) \right)^{1/2} \right]
\end{multline}
where $\ell = f/g$, $L = \log \ell$ and the $Y_i$'s are i.i.d.\ $\sim g$. When $\phi \equiv 1$ and $f = f|_A$ for some measurable set $A \subset \R^d$, then~\eqref{eq:bound-CD} becomes for $n \geq e^{D(f|_A || g)}$
\begin{multline} \label{eq:bound-CD-conditional}
	\E \left( \left \lvert \frac{1}{p_f(A)n} \sum_{i=1}^n \ell(Y_i) \xi_A(Y_i) - 1 \right \rvert \right) \leq \left(\frac{e^{D(f|_A || g)}}{n}\right)^{1/4}\\
 + 2 \left( \P_{f|_A} \left(L_A(X) \geq \frac{1}{2} \log n + \frac{1}{2}D(f|_A || g) \right) \right)^{1/2}
\end{multline}
with $L_A = \log(f|_A / g)$. In the sequel,~\eqref{eq:bound-CD} will be referred to as the CD bound, while~\eqref{eq:bound-CD-conditional}, which as we have just seen is simply a special case of~\eqref{eq:bound-CD}, will be referred to as the conditional CD bound.

An important insight from the CD bounds~\eqref{eq:bound-CD} and~\eqref{eq:bound-CD-conditional} is that in order to show that some auxiliary distribution $g$ is efficient in high dimension for $A$, it is sufficient to control its KL divergence with $f$ and $f|_A$, and also the tails of the log-likelihoods $\log(f/g)$ and $\log(f|_A/g)$ (under $f$ and $f|_A$, respectively). Recall in the next statement that $g$ can be random.

\begin{lemma} \label{lemma:conditions-HDE}
	If $D(f || g)$ is bounded whp and $\P_f(L(X) - D(f||g)\geq t \mid g) \Rightarrow 0$ as $d \to \infty$ for any sequence $t = t(d) \to \infty$, then~\eqref{eq:HDE-WD} holds.

	If $D(f|_A || g)$ is bounded whp and $\P_{f|_A}(L_A(X) - D(f|_A || g)\geq t \mid g) \Rightarrow 0$ as $d \to \infty$ for any sequence $t = t(d) \to \infty$, then~\eqref{eq:HDE-L1} holds.
\end{lemma}

\begin{proof}
 The second part of the lemma follows directly from the conditional CD bound~\eqref{eq:bound-CD-conditional} and with $t = (\log n - D(f|_A||g))/2$ which diverges to $\infty$ since $D(f|_A||g)$ is bounded whp. Note that we can invoke the bound~\eqref{eq:bound-CD-conditional} since $n \geq e^{D(f|_A||g)}$ holds whp, since again $D(f|_A||g)$ is bounded whp. As for the first part of the lemma, the (unconditional) CD bound~\eqref{eq:bound-CD} with $\phi \equiv 1$ implies that $\frac{1}{n} \sum_{i=1}^n \ell(Y_i) \to 1$ in $L_1$ by the same arguments, i.e, $\varepsilon_n \to 0$ with
 \[ \varepsilon_n = \E \left \lvert \frac{1}{n} \sum_{i=1}^n \ell(Y_i) - 1\right \rvert. \]
 According to Theorem 2.3 in~\cite{Chatterjee18:0}, we have
 \[ \E \left ( \frac{\max_{i=1, \ldots, n} \ell(Y_i)}{\sum_{i=1, \ldots, n} \ell(Y_i)} \right) \leq C \max \left \{ \frac{1}{n}, \frac{\log \log (1/\varepsilon_n)}{\log(1/\varepsilon_n)} \right\} \]
 with $C$ a universal constant. Since $\varepsilon_n \to 0$, the bound of the previous vanishes which implies~\eqref{eq:HDE-WD} as desired.
\end{proof}

\begin{rk} \label{rk:WD}
    A simpler condition for~\eqref{eq:HDE-WD}, which does not require to go through the CD bounds, is that there exists $\alpha > 1$ such that $\sup_d \E_g(\ell(X)^\alpha) < \infty$: under this condition, it is easy to prove that $\frac{1}{n} \sum_i \ell(Y_i)$ is tight and that $\frac{1}{n} \max_i \ell(Y_i) \Rightarrow 0$ (where the $Y_i$'s are i.i.d.\ distributed according to $g$), which readily implies~\eqref{eq:HDE-WD}. In Lemma~\ref{lemma:tail-exp} below, we will derive a bound on the $\alpha$-th moment of the likelihood ratio: as this bound also involves the terms $D(f || g)$, $\Sigma$ and $\mu$, going through Lemma~\ref{lemma:tail-exp} rather than the CD bounds does not lead to any significant simplification of the arguments above.
\end{rk}

\subsection{General formula for the Kullback--Leibler divergence}

For the following result, recall that $g|_B = g \xi_B / p_g(B)$ is the measure $g$ conditioned on $B$, and that $\mu^g_B$ and $\Sigma^g_B$ denote the mean and variance of $g|_B$ (see~\eqref{eq:def-mu-sigma-g-B}).

\begin{lemma} \label{lemma:D-g-mid-A-g'}
	Let $g = N(\mu, \Sigma)$ and $g' = N(\mu', \Sigma')$ be two $d$-dimensional Gaussian distributions with $\mu, \mu' \in \R^d$ and $\Sigma, \Sigma' \in \calS_d$, and let $B \subset \R^d$ be any measurable set. Then we have
	\begin{multline} \label{eq:D-g-mid-B-g'-tmp}
		D(g|_B || g') = -\log p_g(B) - \Psi(\Sigma^{-1} \Sigma^g_B) - \frac{1}{2} \left( \mu - \mu^g_B \right)^\top \Sigma^{-1} (\mu - \mu^g_B)\\
		 + \Psi(\Sigma'^{-1} \Sigma^g_B) + \frac{1}{2} \left( \mu' - \mu^g_B \right)^\top \Sigma'^{-1} (\mu' - \mu^g_B).
	\end{multline}
\end{lemma}

\begin{proof}
	By definition, we have
	\[ D(g|_B || g') = \E_{g|_B} \left( \log \left( \frac{g|_B(X)}{g'(X)} \right) \right) = \E_{g|_B} \left( \log \left( \frac{g(X)}{p_g(B) g'(X)} \right) \right) \]
	using for the second equality that the random variable $\xi_B(X)$ is $\P_{g|_B}$-almost surely equal to $1$. Continuing, we get
	\begin{equation} \label{eq:D-g-g'}
		D(g|_B || g') = - \log p_g(B) + \E_{g|_B} \left( \log \left( g(X) \right) \right) - \E_{g|_B} \left( \log \left( g'(X) \right) \right).
	\end{equation}
	We have
	\begin{multline} \label{eq:log-g'}
		\E_{g|_B} \left( \log \left( g'(X) \right) \right) = -\frac{d}{2} \log (2\pi) - \frac{1}{2} \log \det(\Sigma')\\
		- \frac{1}{2} \E_{g|_B} \left( (X-\mu')^\top \Sigma'^{-1}(X - \mu') \right).
	\end{multline}
	Using the identity $\tr(x y^\top) = x^\top y$ and the linearity of the trace and the expectation, which makes them commute, we obtain
	\[ \E_{g|_B} \left( (X-\mu')^\top \Sigma'^{-1}(X - \mu') \right) = \tr \left[ \Sigma'^{-1} \E_{g|_B} \left( (X - \mu') (X-\mu')^\top \right) \right]. \]
	Further, since $\Sigma^g_B$ is the variance of $X$ under $\P_{g|_B}$ and $\mu^g_B$ its mean, we have
	\[ \E_{g|_B} \left( (X - \mu') (X-\mu')^\top \right) = \Sigma^g_B + \left( \mu^g_B -  \mu' \right) \left( \mu^g_B -  \mu' \right)^\top \]
	and so (using again $\tr(V x x^\top) = x^\top V x$)
	\begin{align*}
		\E_{g|_B} \left( (X-\mu')^\top \Sigma'^{-1}(X - \mu') \right) & = \tr \left[ \Sigma'^{-1} \left( \Sigma^g_B + \left( \mu^g_B -  \mu' \right) \left( \mu^g_B -  \mu' \right)^\top \right) \right]\\
		& = \tr \left( \Sigma'^{-1} \Sigma^g_B \right) + (\mu^g_B - \mu')^\top \Sigma'^{-1} (\mu^g_B - \mu').
	\end{align*}
	Plugging in this relation into~\eqref{eq:log-g'}, we obtain
	\begin{multline*}
		\E_{g|_B} \left( \log \left( g'(X) \right) \right) = -\frac{d}{2} \log (2\pi) - \frac{1}{2} \log \det(\Sigma')\\
		- \frac{1}{2} \tr \left( \Sigma'^{-1} \Sigma^g_B \right) - \frac{1}{2} (\mu^g_B - \mu') \Sigma'^{-1} (\mu^g_B - \mu')^\top
	\end{multline*}
	and going back to the definition~\eqref{eq:Psi} of $\Psi$, this gives
	\begin{multline*}
		\E_{g|_B} \left( \log \left( g'(X) \right) \right) = -\frac{d}{2} \log (2\pi) - \Psi(\Sigma'^{-1} \Sigma^g_B) - \frac{1}{2} \log \det(\Sigma^g_B) - \frac{d}{2}\\
		- \frac{1}{2} (\mu^g_B - \mu') \Sigma'^{-1} (\mu^g_B - \mu')^\top.
	\end{multline*}
	Since this formula is valid for any $\mu'$ and $\Sigma'$, it is also valid for $\mu' = \mu$ and $\Sigma' = \Sigma$, and for this choice it gives
	\begin{multline*}
		\E_{g|_B} \left( \log \left( g(X) \right) \right) = -\frac{d}{2} \log (2\pi) - \Psi(\Sigma^{-1} \Sigma^g_B) - \frac{1}{2} \log \det(\Sigma^g_B) - \frac{d}{2}\\
		- \frac{1}{2} (\mu^g_B - \mu) \Sigma^{-1} (\mu^g_B - \mu)^\top.
	\end{multline*}
	Plugging in the two previous relations into~\eqref{eq:D-g-g'} leads to~\eqref{eq:D-g-mid-B-g'-tmp} as desired.
\end{proof}

\begin{corollary}\label{cor:jensen}
	Let $g = N(\mu, \Sigma)$ with $\mu \in \R^d$ and $\Sigma \in \calS_d$. Then for any measurable set $B \in \R^d$, we have
	\[ p_f(B) \geq p_g(B) \exp \left( - \Psi(\Sigma^g_B) - \frac{1}{2} \lVert \mu^g_B \rVert^2 \right). \]
\end{corollary}

\begin{proof}
	We have
	\begin{align*}
		p_f(B) & = \P_f(X \in B)\\
		& = \E_g \left( \frac{f(X)}{g(X)} \xi_B(X) \right)\\
		& = \P_g(X \in B) \E_{g|_B} \left( \frac{f(X)}{g(X)} \right)\\
		& = \E_{g|_B} \left( \frac{f(X)}{g|_B(X)} \right)\\
		& = \E_{g|_B} \left( \exp \left\{ \log \left( \frac{f(X)}{g|_B(X)} \right\} \right) \right).
	\end{align*}
 	Using Jensen's inequality with the convex function $\exp$, we obtain
	\[ p_f(B) \geq \exp \left\{ \E_{g|_B} \left( \log \left( \frac{f(X)}{g|_B(X)} \right) \right) \right\} = \exp \left( - D(g|_B || f) \right). \]
	Applying~\eqref{eq:D-g-mid-B-g'-tmp} with $g' = f$, we see that
	\begin{multline*}
		D(g|_B || f) = -\log p_g(B) - \Psi(\Sigma^{-1} \Sigma^g_B) - \frac{1}{2} \left( \mu - \mu^g_B \right)^\top \Sigma^{-1} (\mu - \mu^g_B)\\
		 + \Psi(\Sigma^g_B) + \frac{1}{2} \lVert \mu^g_B \rVert^2
	\end{multline*}
	and so
	\[ D(g|_B || f) \leq -\log p_g(B) + \Psi(\Sigma^g_B) + \frac{1}{2} \lVert \mu^g_B \rVert^2. \]
	Plugging this inequality in the inequality $p_f(B) \geq e^{-D(g|_B||f)}$ derived above gives the result.
\end{proof}

\subsection{Results on the function $\Psi$}

In this section we gather useful results on $\Psi(\Sigma)$.

\begin{lemma}\label{lemma:bound-Psi-frob}
	There exist two families of positive constants $\{c^\pm_{\varepsilon, K} : \varepsilon, K \in (0,\infty)\}$, independent of the dimension $d$, such that for any $d \geq 1$ and any $\Sigma \in \calS_d$, the following implication holds for any $\varepsilon, K \in (0,\infty)$:
	\[ \varepsilon \leq \lambda_1(\Sigma) \leq \lambda_d(\Sigma) \leq K \Longrightarrow c^-_{\varepsilon, K} \frob{\Sigma - I}^2 \leq \Psi(\Sigma) \leq c^+_{\varepsilon, K} \frob{\Sigma - I}^2. \]
\end{lemma}

\begin{proof}
	Let $\varepsilon, K \in (0,\infty)$ and
 \[ c^-_{\varepsilon, K} = \inf_{\varepsilon \leq x \leq K} \frac{\psi(x)}{(1-x)^2} \ \text{ and } \ c^+_{\varepsilon, K} = \sup_{\varepsilon \leq x \leq K} \frac{\psi(x)}{(1-x)^2}. \]
 Since $\psi(x) \sim \frac{1}{2} (1-x)^2$ for $x \to 1$, $\psi$ is continuous and $\psi(x) > 0$ for any $x > 0$, it follows that $c^-_{\varepsilon, K}$ and $c^+_{\varepsilon, K}$ are strictly positive and finite. Moreover, by definition they satisfy $c^-_{\varepsilon, K} (1-x)^2 \leq \psi(x) \leq c^+_{\varepsilon, K} (1-x)^2$ for any $x \in [\varepsilon, K]$, which gives the result since $\Psi(\Sigma) = \sum_i \psi(\lambda_i(\Sigma))$ and $\frob{\Sigma-I}^2 = \sum_i (\lambda_i(\Sigma) - 1)^2$.
\end{proof}

For the next statement, recall that a sequence of real-valued random variables $X$ is said to be bounded whp if $\P(\lvert X \rvert \leq K) \to 1$ for some $K \geq 0$.

\begin{lemma} \label{lemma:equivalences}
 For each $d$ consider $\Sigma \in \calS_d$ possibly random. Then the following three conditions are equivalent:
 \begin{enumerate}
 \item $\Psi(\Sigma)$ is bounded whp;
 \item $\Psi(\Sigma^{-1})$ is bounded whp;
 \item the three sequences $1/\lambda_1(\Sigma)$, $\lambda_d(\Sigma)$ and $\frob{\Sigma - I}$ are bounded whp.
 \end{enumerate}
\end{lemma}

\begin{proof}
	Let us first prove these equivalences using the notion of almost-sure boundedness instead of the notion of boundedness whp: we will relax this assumption at the end of the proof. Let us first show that $1 \Rightarrow 3$ assuming that $\Psi(\Sigma)$ is almost surely bounded, i.e., assume that $\Psi(\Sigma)$ is almost surely bounded and let us show that $1/\lambda_1(\Sigma)$, $\lambda_d(\Sigma)$ and $\frob{\Sigma - I}$ are almost surely bounded. We have $\Psi(\Sigma) \geq \psi(\lambda_1(\Sigma))$ and so $\sup_d \psi(\lambda_1(\Sigma)) < \infty$, and so necessarily $\inf_d \lambda_1(\Sigma) > 0$ because $\psi(x) \to \infty$ as $x \to 0$. The same argument implies $\sup_d \lambda_d(\Sigma) < \infty$. And since $\lambda_1(\Sigma)$ and $\lambda_d(\Sigma)$ are bounded away from $0$ and $\infty$, the boundedness of $\frob{\Sigma - I}$ comes from Lemma~\ref{lemma:bound-Psi-frob}.
	
	The implication $3 \Rightarrow 1$ is immediate in view of Lemma~\ref{lemma:bound-Psi-frob}.
	
	To conclude, note that since $\Sigma^{-1} - I$ is symmetric, the square of its Frobenius norm is equal to the sum of its eigenvalues:
	\begin{equation} \label{eq:bound--1}
		\frob{\Sigma^{-1} - I}^2 = \sum_i \left( \frac{1}{\lambda_i(\Sigma)} - 1 \right)^2 \leq \frac{1}{\lambda_1(\Sigma)^2} \frob{\Sigma - I}^2.
	\end{equation}
	In particular, $3$ implies that $\frob{\Sigma^{-1} - I}$ is almost surely bounded and so we can invoke the implication $3 \Rightarrow 1$ for $\Sigma^{-1}$, which shows that $3 \Rightarrow 2$. The implication $2 \Rightarrow 3$ follows for similar reasons.
	
	Let us now show that our results apply assuming that the random variables are bounded whp. Let us for instance show that $1 \Rightarrow 2$, the arguments for the other implications are the same. Let $K \geq 0$ such that $\P(\Psi(\Sigma) \leq K) \to 1$. Then under $\P(\cdot \mid \Psi(\Sigma) \leq K)$, $\Psi(\Sigma)$ is almost surely bounded (by $K$) and so we have proved that $\Psi(\Sigma^{-1})$ is almost surely bounded, i.e., there exists $K'$ such that $\P(\Psi(\Sigma^{-1}) \leq K' \mid \Psi(\Sigma) \leq K) = 1$. Writing
	\begin{multline*}
		\P(\Psi(\Sigma^{-1}) \leq K') = \P(\Psi(\Sigma) \leq K) \P(\Psi(\Sigma^{-1}) \leq K' \mid \Psi(\Sigma) \leq K)\\
		+ \P(\Psi(\Sigma^{-1}) \leq K', \Psi(\Sigma) > K)
	\end{multline*}
	and noting that $\P(\Psi(\Sigma) \leq K) \to 1$, $\P(\Psi(\Sigma^{-1}) \leq K' \mid \Psi(\Sigma) \leq K) = 1$ and $\P(\Psi(\Sigma^{-1}) \leq K', \Psi(\Sigma) > K) \to 0$, we obtain that $\P(\Psi(\Sigma^{-1}) \leq K') \to 1$ as desired.
\end{proof}

\begin{lemma} \label{lemma:boundedness-Sigma-Sigma'}
	For each $d \geq 1$, let $\Sigma, \Sigma' \in \calS_d$ possibly random. If $\Psi(\Sigma)$ and $\Psi(\Sigma')$ are bounded whp, then $\Psi(\Sigma \Sigma')$ is bounded whp.
\end{lemma}

\begin{proof}
	As in Lemma~\ref{lemma:equivalences}, we prove the results with the notion of almost-sure boundedness instead of boundedness whp. For simplicity, the almost surely quantifiers are left out. So assume that $\Psi(\Sigma)$ and $\Psi(\Sigma')$ are bounded, and let us show that $\Psi(\Sigma \Sigma')$ is bounded. Lemma~\ref{lemma:equivalences} implies that the sequences $\lambda_d(\Sigma)$, $1/\lambda_1(\Sigma)$ and $\frob{\Sigma - I}$ are bounded, and the same holds with $\Sigma'$ instead of $\Sigma$. Since $\lambda_d(\Sigma)$ is the matrix-norm induced by the $L_2$-norm on $\R^d$, it is submultiplicative, and so $\lambda_d(\Sigma \Sigma') \leq \lambda_d(\Sigma) \lambda_d(\Sigma')$, which implies $\lambda_1(\Sigma \Sigma') \geq \lambda_1(\Sigma) \lambda_1(\Sigma')$ since $\lambda_1(\Sigma) = 1/\lambda_d(\Sigma^{-1})$. Therefore, $\lambda_d(\Sigma \Sigma')$ and $1/\lambda_1(\Sigma \Sigma')$ are bounded. Moreover, since
	\begin{align*}
		\frob{\Sigma \Sigma' - I} = \frob{(\Sigma - I)(\Sigma' - I) + \Sigma' - I + \Sigma - I},
	\end{align*}
	the triangle inequality and the sub-multiplicativity of the Frobenius norm imply that
	\[ \frob{\Sigma \Sigma' - I} \leq \frob{\Sigma - I} \frob{\Sigma' - I} + \frob{\Sigma - I} + \frob{\Sigma' - I} \]
	and so $\frob{\Sigma \Sigma' - I}$ is bounded. Lemma~\ref{lemma:equivalences} implies that $\Psi(\Sigma \Sigma')$ is bounded.
\end{proof}

\begin{corollary} \label{cor:boundedness-sigma-mu}
	Let $g = N(\mu, \Sigma)$ with $\mu \in \R^d$ and $\Sigma \in \calS_d$, $B \subset \R^d$ measurable and $g_B = N(\mu^g_B, \Sigma^g_B)$; $\mu$, $\Sigma$ and $B$ may be random. If $\lVert \mu \rVert$, $\Psi(\Sigma)$ and $1/p_g(B)$ are bounded whp, then $D(g|_B || g_B)$, $\Psi(\Sigma^g_B)$ and $\lVert \mu^g_B \rVert$ are bounded whp.
\end{corollary}

\begin{proof}
	As in the previous two proofs, we prove the results with the notion of almost-sure boundedness instead of boundedness whp. If we apply~\eqref{eq:D-g-mid-B-g'-tmp} with $\mu' = \mu^g_B$ and $\Sigma' = \Sigma^g_B$, we obtain
	\[ D(g|_B || g_B) = -\log p_g(B) - \Psi(\Sigma^{-1} \Sigma^g_B) - \frac{1}{2} \left( \mu - \mu^g_B \right)^\top \Sigma^{-1} (\mu - \mu^g_B). \]
	Since $\Psi \geq 0$, we see that $D(g|_B || g_B) \leq -\log p_g(B)$ which gives the boundedness of $D(g|_B || g_B)$. Moreover, since $D(g|_B || g_B) \geq 0$ we obtain
	\[ \Psi(\Sigma^{-1} \Sigma^g_B) + \frac{1}{2} \left( \mu - \mu^g_B \right)^\top \Sigma^{-1} (\mu - \mu^g_B) \leq -\log p_g(B). \]
	This shows that $\Psi(\Sigma^{-1} \Sigma^g_B)$ is bounded. But $\Psi(\Sigma)$ is assumed to be bounded, and so $\Psi(\Sigma^{-1})$ is bounded by Lemma~\ref{lemma:equivalences}, which implies that $\Psi(\Sigma^g_B)$ is bounded by Lemma~\ref{lemma:boundedness-Sigma-Sigma'}. Likewise, the boundedness of $\left( \mu - \mu^g_B \right)^\top \Sigma^{-1} (\mu - \mu^g_B)$ implies that of $\lVert \mu^g_B \rVert$ because:
	\[ \left( \mu - \mu^g_B \right)^\top \Sigma^{-1} (\mu - \mu^g_B) \geq \frac{1}{\lambda_d(\Sigma)} \lVert \mu - \mu^g_B \rVert^2 \geq \frac{1}{\lambda_d(\Sigma)} \left( \lVert \mu^g_B \rVert - \lVert \mu \rVert \right)^2. \]
	Since $\lambda_d(\Sigma)$ is bounded (by Lemma~\ref{lemma:equivalences}, because $\Psi(\Sigma)$ is), and $\lVert \mu \rVert$ is bounded by assumption, the boundedness of $\left( \mu - \mu^g_B \right)^\top \Sigma^{-1} (\mu - \mu^g_B)$ indeed implies that of $\lVert \mu^g_B \rVert$ by the inequality of the previous display.
\end{proof}

\subsection{Bound on the tail of the log-likelihoods}

In the next statement recall that $f = N(0,I)$ is the standard Gaussian density in dimension $d$ and that $f|_B = f \xi_B / p_f(B)$ is the density $f$ conditioned on $B$ with mean $\mu_B$ and variance $\Sigma_B$ (see~\eqref{eq:def-mu-sigma-f-B}).

\begin{lemma} \label{lemma:bound-var}
	For $B \subset \R^d$ measurable, $y \in \R^d$ and $V \in \calM_d$ symmetric, we have
	\[ \Var_{f|_B}(y^\top X) \leq \lambda_d(\Sigma_B) \lVert y \rVert^2 \ \text{ and } \ \Var_{f|_B} \left( X^\top V X \right) \leq \frac{2}{p_f(B)} \frob{V}^2. \]
\end{lemma}

\begin{proof}
	The first inequality follows from the fact that $\Var_{f|_B}(y^\top X) = y^\top \Sigma_B y$ and the variational characterization of eigenvalues. Let us prove the second inequality. First, note that for any measurable function $h: \R^d \to \R$, we have $\Var_{f|_B}(h(X)) \leq \frac{1}{p_f(B)} \Var_f(h(X))$. Indeed, we have
	\[ \Var_{f|_B}(h(X)) = \E_{f|_B} \left[ (h(X) - \E_{f|_B}(h(X))^2 \right] \leq \E_{f|_B} \left[ (h(X) - \E_f(h(X))^2 \right] \]
	where the last inequality follows by the variational characterization of the mean. By definition of $f|_B$, we have
	\[ \E_{f|_B} \left[ (h(X) - \E_f(h(X)))^2 \right] = \E_f \left[ (h(X) - \E_f(h(X)))^2 \mid X \in B \right] \]
	from which the desired inequality $\Var_{f|_B}(h(X)) \leq \frac{1}{p_f(B)} \Var_f(h(X))$ readily follows. In particular,
	\[ \Var_{f|_B} \left( X^\top V X \right) \leq \frac{1}{p_f(B)} \Var_f \left( X^\top V X \right). \]
	 Write $V = U^\top \Delta U$ with $U$ orthonormal and $\Delta$ the diagonal matrix with diagonal elements the $\lambda_i(V)$'s, so that $X^\top V X = (UX)^\top \Delta (UX)$. Under $\P_f$, $X$ is standard Gaussian and since $U$ is orthonormal, $UX$ is also standard Gaussian, so that
	 \[ \Var_f( X^\top V X ) = \Var_f(X^\top \Delta X). \]
	 Since $X^\top \Delta X = \sum_i \lambda_i(V) X(i)^2$ with, under $\P_f$, the $X(i)$'s i.i.d., we obtain
	 \[ \Var_f( X^\top V X ) = \Var_f(X(1)^2) \sum_i \lambda_i(V)^2 = 2 \frob{V}^2 \]
	using for the last equality that $\sum_i \lambda_i(V)^2 = \frob{V}^2$ and $\Var_f(X(1)^2) = 2$. This proves the result.
\end{proof}

\begin{corollary} \label{cor:quadratic-bound-tail}
	Let $B \subset \R^d$ measurable, $g = N(\mu, \Sigma)$ with $\mu \in \R^d$ and $\Sigma \in \calS_d$ and $L = \log(f|_B / g)$. Then for any $t > 0$ we have
	\begin{equation} \label{eq:tail}
		\P_{f|_B}(L(X) - \E_{f|_B}(L(X)) \geq t) \leq \frac{4}{t^2} \left( \frac{2\frob{\Sigma^{-1} - I}^2}{p_f(B)} + \frac{\lambda_d(\Sigma_B)}{\lambda_1(\Sigma)^2} \lVert \mu \rVert^2 \right).
	\end{equation}
\end{corollary}

\begin{proof}
	For $x \in B$, we have
	\begin{align*}
		L(x) & = \log (f|_B(x) / g(x))\\
		& = - \log p_f(B) - \frac{1}{2} \lVert x \rVert^2 + \frac{1}{2} \log \det(\Sigma) + \frac{1}{2} (x - \mu)^\top \Sigma^{-1} (x - \mu)\\
		& = - \log p_f(B) + \frac{1}{2} \mu^\top \Sigma^{-1} \mu + \frac{1}{2} \log \det(\Sigma) + \frac{1}{2} x^\top (\Sigma^{-1} - I) x - x^\top \Sigma^{-1} \mu.
	\end{align*}
	Let $Z_1 = \frac{1}{2} X^\top (\Sigma^{-1} - I) X$ and $Z_2 = -X^\top \Sigma^{-1} \mu$, and $\bar Z_i = Z_i - \E_{f|_B}(Z_i)$ for $i = 1, 2$ be their centered versions: then $L(X) - \E_{f|_B}(L(X)) = \bar Z_1 + \bar Z_2$ and so
	\begin{align*}
		\P_{f|_B}(L(X) - \E_{f|_B}(L(X)) \geq t) & = \P_{f|_B}( \bar Z_1 + \bar Z_2 \geq t)\\
		& \leq \P_{f|_B}( \bar Z_1 \geq t/2) + \P( \bar Z_2 \geq t/2)\\
		& \leq \frac{4}{t^2} \left( \Var_{f|_B}(Z_1) + \Var_{f|_B}(Z_2) \right)
	\end{align*}
	and so Lemma~\ref{lemma:bound-var} gives
	\[ \P_{f|_B}(L(X) - \E_{f|_B}(L(X)) \geq t) \leq \frac{4}{t^2} \left( \frac{2}{p_f(B)} \frob{\Sigma^{-1} - I}^2 + \lambda_d(\Sigma_B) \lVert \Sigma^{-1} \mu \rVert^2 \right). \]
	The result thus follows from the fact that $\lVert \Sigma^{-1} \mu \rVert^2 = \mu^\top \Sigma^{-2} \mu \leq \lambda_d(\Sigma^{-2}) \lVert \mu \rVert^2$ and $\lambda_d(\Sigma^{-2}) = 1/\lambda_1(\Sigma)^2$.
\end{proof}

When $B = \R^d$, we will sometimes need the following strengthening of Corollary~\ref{cor:quadratic-bound-tail}. In the sequel, let $\alpha_*(\Sigma)$ for $\Sigma \in \calS_d$ be defined as follows:
\begin{equation} \label{eq:alpha*}
	\alpha_*(\Sigma) = \min \left(1, \frac{\lambda_1(\Sigma)}{1-\lambda_1(\Sigma)} \right) = \left\{ \begin{array}{cl}
 \displaystyle \frac{\lambda_1(\Sigma)}{1-\lambda_1(\Sigma)} & \text{ if } \lambda_1(\Sigma) < \frac{1}{2}, \\
 1 & \text{ else}
\end{array} \right.
\end{equation}

\begin{lemma}\label{lem:alpha*}
	If $\Sigma \in \calS_d$ and $\alpha < \alpha_*(\Sigma)$, then $(\alpha + 1) I - \alpha \Sigma^{-1} \in \calS_d$.
\end{lemma}

\begin{proof}
	Let $W = (\alpha + 1) I - \alpha \Sigma^{-1}$: by definition, it is symmetric and so we only have to show that $\lambda_1(W) > 0$. We have
    \[ \lambda_1(W) = \alpha + 1 + \lambda_1(-\alpha \Sigma^{-1}) = \alpha + 1 - \alpha \lambda_d(\Sigma^{-1}) = \alpha + 1 - \frac{\alpha}{\lambda_1(\Sigma)} \]
	and so
    \[ \lambda_1(W) = 1 + \frac{\lambda_1(\Sigma)-1}{\lambda_1(\Sigma)} \alpha = \frac{1-\lambda_1(\Sigma)}{\lambda_1(\Sigma)} \left( \frac{\lambda_1(\Sigma)}{1-\lambda_1(\Sigma)} - \alpha \right). \]
	The first equality clearly shows that $\lambda_1(W) > 0$ if $\lambda_1(\Sigma) \geq 1$. For $\lambda_1(\Sigma) < 1/2$, the second equality can be rewritten as $\lambda_1(W) = (\alpha_*(\Sigma) - \alpha) / \alpha_*(\Sigma)$ which is $>0$. Finally, for $\lambda_1(\Sigma) \in [1/2,1)$, we have $\frac{\lambda_1(\Sigma)}{1-\lambda_1(\Sigma)} \geq 1 = \alpha_*(\Sigma)$ and so using that $(1-\lambda_1(\Sigma))/\lambda_1(\Sigma) > 0$, the second inequality leads to
    \[ \lambda_1(W) \geq \frac{1-\lambda_1(\Sigma)}{\lambda_1(\Sigma)}(1-\alpha) = \frac{1-\lambda_1(\Sigma)}{\lambda_1(\Sigma)}(\alpha_*(\Sigma)-\alpha) > 0. \]
	This proves the result.
\end{proof}

\begin{lemma} \label{lemma:tail-exp}
	Let $g = N(\mu, \Sigma)$ with $\mu \in \R^d$ and $\Sigma \in \calS_d$ and $L = \log(f/g)$. Then for every $\alpha < \alpha' < \alpha_*(\Sigma)$, we have
	\begin{multline} \label{eq:tail-exp}
		\E_f \left[ \left( \frac{f(X)}{g(X)} \right)^\alpha \right]\\
		\leq \exp \left(\alpha D(f || g) + \frac{1}{2} q \alpha^2 \lVert \Sigma^{-1} \mu \rVert^2 + \frac{\alpha}{2 \alpha'} \Psi((\alpha'+1)I - \alpha' \Sigma^{-1}) \right)
	\end{multline}
	where $q = \alpha'/(\alpha'-\alpha)$.
\end{lemma}

\begin{proof}
	Let $W = (\alpha'+1) I - \alpha' \Sigma^{-1}$, which belongs to $\calS_d$ by Lemma~\ref{lem:alpha*} (so that $\Psi(W)$ is well defined). Let $\bar Z_1 = \frac{1}{2} X^\top (\Sigma^{-1} - I) X - \frac{1}{2} \tr(\Sigma^{-1} - I)$ and $\bar Z_2 = - X^\top \Sigma^{-1} \mu$: proceeding similarly as in the proof of Corollary~\ref{cor:quadratic-bound-tail}, we see that $L(X) - \E_f(L(X)) = \bar Z_1 + \bar Z_2$ and so
	\begin{align*}
		\E_f \left[ \left( \frac{f(X)}{g(X)} \right)^\alpha \right] & = \E_f \left[ \exp \left( \alpha L(X) \right) \right]\\
		 & = e^{\alpha D(f || g)} \E_f \left[ \exp \left( \alpha (L(X) - D(f || g)) \right) \right]\\
		 & = e^{\alpha D(f || g)} \E_f \left( e^{\alpha \bar Z_1} e^{\alpha \bar Z_2} \right).
	\end{align*}
	Let $p = \alpha'/\alpha$ and $q = p/(p-1) = \alpha'/(\alpha'-\alpha)$: then $1/p + 1/q = 1$ and so H\"older's inequality gives
	\[ \E_f \left[ \left( \frac{f(X)}{g(X)} \right)^\alpha \right] \leq e^{\alpha D(f || g)} \left\{ \E_f \left( e^{p \alpha \bar Z_1} \right) \right\}^{1/p} \left\{ \E_f \left( e^{q \alpha \bar Z_2} \right) \right\}^{1/q}. \]
	Recall that $\bar Z_2 = - X^\top \Sigma^{-1} \mu$: since $\E_f(e^{x^\top X}) = e^{\frac{1}{2} \lVert x \rVert^2}$ for any $x \in \R^d$, we obtain
	\[ \left\{ \E_f(e^{q\alpha \bar Z_2}) \right\}^{1/q} = \left\{ \E_f(e^{-q\alpha \mu^\top \Sigma^{-1}X}) \right\}^{1/q} = e^{\frac{1}{2q} \lVert q\alpha \mu^\top \Sigma^{-1} \rVert^2} = e^{\frac{1}{2} q \alpha^2 \mu^\top \Sigma^{-2} \mu}. \]
	Let us now control the exponential moment of $\bar Z_1$. We have
	\begin{align*}
	 \E_f(e^{p \alpha \bar Z_1}) & = \E_f(e^{\alpha' \bar Z_1})\\
	 & = e^{- \frac{1}{2} \alpha' \tr(\Sigma^{-1}-I)} \E_f(e^{\frac{1}{2} \alpha' X^\top (\Sigma^{-1}-I) X})\\
	 & = e^{\frac{1}{2} \tr(W-I)} \int \frac{1}{(2 \pi)^{d/2}} e^{\frac{1}{2} \alpha' x^\top (\Sigma^{-1}-I) x - \frac{1}{2} x^\top x} \d x\\
	 & = e^{\frac{1}{2} \tr(W-I)} \int \frac{1}{(2 \pi)^{d/2}} e^{-\frac{1}{2} x^\top W x} \d x.
	\end{align*}
	Since we have seen that $W \in \calS_d$, we have
	\[ \int \frac{\det(W)^{1/2}}{(2 \pi)^{d/2}} e^{-\frac{1}{2} x^\top W x} \d x = 1 \]
	and so
	\[ \left\{ \E(e^{p \alpha \bar Z_1}) \right\}^{1/p} = \exp \left( \frac{1}{2p} \tr(W-I) - \frac{1}{2p} \log \det(W) \right) = e^{\frac{1}{p} \Psi(W)}. \]
	Gathering the previous bounds leads to the desired result.
\end{proof}

\subsection{A sufficient condition for high-dimensional efficiency}

The following result identifies conditions under which~\eqref{eq:HDE-WD} and~\eqref{eq:HDE-L1} hold for a Gaussian density $g = N(\mu, \Sigma)$. It shows in particular that~\eqref{eq:HDE-L1} is slightly more demanding than~\eqref{eq:HDE-WD}: for~\eqref{eq:HDE-WD}, it is enough that $\Psi(\Sigma)$ and $\lVert \mu \rVert$ are bounded whp (note in particular that this condition does not depend on $A$), and for~\eqref{eq:HDE-L1}, one needs in addition that $1/p_f(A)$ is bounded.

An intuitive interpretation of these conditions is as follows. Since
\begin{equation} \label{eq:D-f-g}
	D(f || g) = \Psi(\Sigma^{-1}) + \frac{1}{2} \mu^\top \Sigma^{-1} \mu,
\end{equation}
the assumption that $\Psi(\Sigma)$ and $\lVert \mu \rVert$ are bounded means that $g$ remains close to $f$. On the other hand, since $D(f|_A||f) = - \log p_f(A)$, the assumption $1/p_f(A)$ bounded means that $f|_A$ remains close to $f$.

\begin{prop} \label{prop:first-step}
	Let $\mu \in \R^d$, $\Sigma \in \calS_d$ and $B \subset \R^d$ measurable ($\mu$, $\Sigma$ and $B$ may be random) and $g = N(\mu, \Sigma)$. Then the following holds:
	\begin{itemize}
		\item if $\Psi(\Sigma)$ and $\lVert \mu \rVert$ are bounded whp, then~\eqref{eq:HDE-WD} holds;
		\item if $\Psi(\Sigma)$, $\lVert \mu \rVert$ and $1/p_f(B)$ are bounded whp, then~\eqref{eq:HDE-L1} holds.
	\end{itemize}
	In particular, if $\Psi(\Sigma)$, $\lVert \mu \rVert$ and $1/p_f(B)$ are bounded whp, then $g = N(\mu, \Sigma)$ is efficient in high dimension for $B$.
\end{prop}

\begin{proof}
	As before, it is enough to prove the result for deterministic $\mu$, $\Sigma$ and $B$, and by assuming that the quantities of interest are bounded. So assume in the rest of the proof that $\Psi(\Sigma)$ and $\lVert \mu \rVert$ are bounded: we first prove that~\eqref{eq:HDE-WD} holds, and then that ~\eqref{eq:HDE-L1} holds under the additional assumption that $1/p_f(B)$ is bounded. The boundedness of $\Psi(\Sigma)$ implies by Lemma~\ref{lemma:equivalences} that $1/\lambda_1(\Sigma)$, $\lambda_d(\Sigma)$, $\Psi(\Sigma^{-1})$ and $\frob{\Sigma^{-1} - I}$ are bounded, which will be used without further notice in the rest of proof. Recall that $L = \log(f/g)$ and that $L_B = \log(f|_B / g)$, with respective means $\E_f(L(X)) = D(f||g)$ and $\E_{f|_B}(L_B(X)) = D(f|_B||g)$.\\
	
	\noindent \textit{Proof of~\eqref{eq:HDE-WD}.} According to Lemma~\ref{lemma:conditions-HDE} it is enough to prove that $D(f||g)$ is bounded and that $\P_f(L(X) - \E_f(L(X)) \geq t) \to 0$ for any sequence $t \to \infty$. Since $\mu^\top \Sigma^{-1} \mu \leq \lVert \mu \rVert^2 / \lambda_1(\Sigma)$, it follows from~\eqref{eq:D-f-g} that $D(f || g)$ is bounded. Let us now control the tail of $L$. Using~\eqref{eq:tail} with $B = \R^d$, we get
	\[ \P_f(L(X) - \E_f(L(X)) \geq t) \leq \frac{4}{t^2} \left( 2 \frob{\Sigma^{-1} - I}^2 + \frac{1}{\lambda_1(\Sigma)^2} \lVert \mu \rVert^2 \right). \]
	The upper bound is thus of the form $C/t$ with $\sup_d C < \infty$, which implies as desired that $\P_f(L(X) - \E_f(L(X)) \geq t) \to 0$ as $d \to \infty$ for any sequence $t \to \infty$.\\

	\noindent \textit{Proof of~\eqref{eq:HDE-L1}.} According to Lemma~\ref{lemma:conditions-HDE} it is enough to prove that $D(f|_B|||g)$ is bounded and that $\P_{f|_B}(L_B(X) - \E_{f|_B}(L_B(X)) \geq t) \to 0$ for any sequence $t \to \infty$. If we apply~\eqref{eq:D-g-mid-B-g'-tmp} with $\mu' = \mu^g_B$ and $\Sigma' = \Sigma^g_B$, we obtain
 	\[ D(g|_B || g_B) = -\log p_g(B) - \Psi(\Sigma^{-1} \Sigma^g_B) - \frac{1}{2} \left( \mu - \mu^g_B \right)^\top \Sigma^{-1} (\mu - \mu^g_B) \]
	and so for any $g'=N(\mu',\Sigma')$, \eqref{eq:D-g-mid-B-g'-tmp} can be rewritten as
	\begin{equation} \label{eq:D-g-mid-B-g'}
		D(g|_B || g') = D(g|_B || g_B) + \Psi(\Sigma'^{-1} \Sigma^g_B) + \frac{1}{2} \left( \mu' - \mu^g_B \right)^\top \Sigma'^{-1} (\mu' - \mu^g_B).
	\end{equation} 
    Plugging $g = f$ and $g' = g$ in this relation, we get
	\[ D(f|_B || g) = D(f|_B || f_B) + \Psi(\Sigma^{-1} \Sigma_B) + \frac{1}{2} \left( \mu - \mu_B \right)^\top \Sigma^{-1} (\mu - \mu_B). \]
	By Corollary~\ref{cor:boundedness-sigma-mu} (with $g = f$, needing $\inf_d p_f(B) > 0$), we see that $D(f|_B || f_B)$, $\Psi(\Sigma_B)$ and $\lVert \mu_B \rVert$ are bounded. Combining the results from Lemmas~\ref{lemma:equivalences} and~\ref{lemma:boundedness-Sigma-Sigma'}, this implies the boundedness of $\Psi(\Sigma^{-1} \Sigma_B)$ and of $\left( \mu - \mu_B \right)^\top \Sigma^{-1} (\mu - \mu_B)$ which proves that $D(f|_B || g)$ is bounded. Let us now turn to controlling the tail of $L_B$. Using~\eqref{eq:tail}, we get
	\[ \P_{f|_B}(L_B(X) - \E_{f|_B}(L_B(X)) \geq t) \leq \frac{4}{t^2} \left( \frac{2 \frob{\Sigma^{-1} - I}^2}{p_f(B)} + \frac{\lambda_d(\Sigma_B)}{\lambda_1(\Sigma)^2} \lVert \mu \rVert^2 \right) \]
	which implies as above that $\P_{f|_B}(L_B(X) - \E_{f|_B}(L_B(X)) \geq t) \to 0$ as $d \to \infty$ for any sequence $t = t(d) \to \infty$. This concludes the proof of the lemma.
\end{proof}

\subsection{Quantiles}

Let us finally mention a last result which will be needed to study the CE scheme. Recall the assumption in Theorem~\ref{thm:CE} that $\varphi: \R^d \to \R$ has no atom, i.e., for every $x \in \R$ the set $\varphi^{-1}(\{x\}) \subset \R^d$ has zero Lebesgue measure.

\begin{lemma} \label{lemma:no-atom}
    Let $\varphi: \R^d \to \R$ measurable, $g$ a $d$-dimensional Gaussian distribution and $F(x) = \P_g(\varphi(X) \leq x)$. If $\varphi$ has no atom, then $F$ is continuous and $F(F^{-1}(x)) = x$ for every $x \in (0,1)$.
\end{lemma}

\begin{proof}
    We have $F(x) - F(x-) = \P_g(\varphi(X) = x) = \P_g(\varphi^{-1}(\{x\}) = 0$ by assumption on~$\varphi$ (and since $g$ is absolutely continuous with respect to Lebesgue measure). The continuity of~$F$ then implies the relation $F(F^{-1}(x)) = x$, see for instance~\cite[Lemma $13.6.4$, Equation $(6.6)$]{Whitt02:0}.
\end{proof}

\section{Proof of Theorem~\ref{thm:g_proj}} \label{sec:proof-g_proj}

\subsection{High-dimensional efficiency of $g_\proj$}

In the rest of this section, we fix the notation as in the statement of Theorem~\ref{thm:g_proj}. According to Lemma~\ref{prop:first-step}, it is enough to prove that $\lVert \mu_A \rVert$ and $\Psi(\Sigma_\proj)$ are bounded. The following lemma will be needed in order to control $\Psi(\Sigma_\proj)$.

\begin{lemma} \label{lemma:proj-*}
	Let $\Sigma \in \calS_d$, $r \leq d$, $(d_k, k = 1, \dots, r)$ an orthonormal family and $\Sigma' \in \calS_d$ defined by
	\[ \Sigma' = \sum_{k=1}^r (v_k - 1) d_k d_k^\top + I \ \text{ with } \ v_k = d_k^\top \Sigma d_k. \]
	Then we have
	\[ \lambda_1(\Sigma') \geq \min(1, \lambda_1(\Sigma)), \ \lambda_d(\Sigma') \leq \max(1, \lambda_d(\Sigma)) \]
	and $\frob{\Sigma' - I} \leq \frob{\Sigma - I}$.
\end{lemma}

\begin{proof}
	Complete the $(d_k, k = 1, \dots, r)$ into an orthonormal basis $(d_k, k = 1, \ldots, d)$. By construction, the eigenvalues of $\Sigma'$ are the $v_k$'s (associated to the $d_k$ for $k = 1, \ldots, r$) and $1$ (associated to the $d_k$ for $k = r+1, \ldots, d$). For any $x \in \R^d$ with $\lVert x \rVert = 1$, we have
	\[ \lambda_1(\Sigma) \leq x^\top \Sigma x \leq \lambda_d(\Sigma) \]
	and since $v_k = d_k^\top \Sigma d_k$ for $k = 1, \ldots, r$, this gives $\lambda_1(\Sigma) \leq v_k \leq \lambda_d(\Sigma)$ for $k = 1, \ldots, r$. Let us show the inequality $\lambda_1(\Sigma') \geq \min(1, \lambda_1(\Sigma))$ by distinguishing two cases:
	\begin{description}
		\item[Case $1$:] if all the $v_k$'s are $\geq 1$, then $\lambda_1(\Sigma') = 1$ and so $\lambda_1(\Sigma') = 1 \geq \min(1, \lambda_1(\Sigma))$ as desired;
		\item[Case $2$:] otherwise, there is some $v_k < 1$, in which case $\lambda_1(\Sigma') = v_i$ for some $i$. But since $v_i \geq \lambda_1(\Sigma)$, the inequality $\lambda_1(\Sigma') \geq \min(1, \lambda_1(\Sigma))$ is also satisfied in this case.
	\end{description}
	A similar discussion shows that $\lambda_d(\Sigma') \leq \max(1, \lambda_d(\Sigma))$. Let us now show that $\frob{\Sigma' - I} \leq \frob{\Sigma - I}$. Since the eigenvalues of $\Sigma'$ are the $v_k$'s and $1$, we have
	\[ \frob{\Sigma' - I}^2 = \sum_i (\lambda_i(\Sigma') - 1)^2 = \sum_{k=1}^r (v_k-1)^2. \]
	By definition of $v_k$,
	\begin{align*}
		\sum_{k=1}^r (v_k-1)^2 & = \sum_{k=1}^r (d_k^\top \Sigma d_k-1)^2\\
		& = \sum_{k=1}^r (d_k^\top (\Sigma - I) d_k)^2\\
		& \leq \sum_{k=1}^d (d_k^\top (\Sigma - I) d_k)^2.
	\end{align*}
	Let $U$ orthonormal such that $\Sigma = U^\top \Lambda U$ with $\Lambda$ the diagonal matrix with diagonal elements the $\lambda_i(\Sigma)$'s. Then $d_k^\top (\Sigma - I) d_k = \tilde d_k^\top (\Lambda - I) \tilde d_k$ with $\tilde d_k = U d_k$. We then have
	\begin{align*}
		\sum_{k=1}^r (v_k-1)^2 & \leq \sum_{k=1}^d (\tilde d_k^\top (\Lambda - I) \tilde d_k)^2\\
		& = \sum_{k=1}^d \left( \sum_{i=1}^d (\tilde d_k(i))^2 \lambda_i(\Sigma - I) \right)^2\\
		& \leq \sum_{k=1}^d \left( \sum_{i=1}^d (\tilde d_k(i))^2 \right) \left( \sum_{i=1}^d (\tilde d_k(i))^2 \lambda_i(\Sigma - I)^2 \right)
	\end{align*}
	using Cauchy--Schwarz for the last inequality (with $\tilde d_k(i)$ on the one hand, and $\tilde d_k(i) \lambda_i(\Sigma - I)$ on the other hand). Since $U$ is orthonormal and the $d_k$'s form an orthonormal basis, the $\tilde d_k$ also form an orthonormal basis, in particular $\sum_{i=1}^d (\tilde d_k(i))^2 = 1$ and so continuing the previous derivation leads to
	\begin{align*}
		\sum_{k=1}^r (v_k-1)^2 & \leq \sum_{k=1}^d \left( \sum_{i=1}^d (\tilde d_k(i))^2 \lambda_i(\Sigma - I)^2 \right)\\
		& = \sum_{i=1}^d \lambda_i(\Sigma - I)^2 \sum_{k=1}^d (\tilde d_k(i))^2\\
		& = \sum_{i=1}^d \lambda_i(\Sigma - I)^2
	\end{align*}
	using $\sum_{k=1}^d (\tilde d_k(i))^2 = 1$ to derive the last equality, which holds because the $\tilde d_k$'s form an orthonormal basis. Since this last quantity is equal to $\frob{\Sigma - I}^2$, this gives the result.
\end{proof}

We get the following corollary, whose first part proves the part of Theorem~\ref{thm:g_proj} related to $g_\proj$.

\begin{corollary} \label{cor:Sigma-proj-mu-A-bounded}
	If $1/p_f(A)$ is bounded, then $g_\proj$ is efficient in high dimension for $A$.
	
	More precisely, if $1/p_f(A)$ is bounded, then $\lVert \mu_A \rVert$, $\Psi(\Sigma_A)$, $\Psi(\Sigma_\proj)$, $1/\lambda_1(\Sigma_\proj)$, $\frob{\Sigma_\proj - I}$ and $\lambda_d(\Sigma_\proj)$ are bounded.
\end{corollary}

\begin{proof}
	The boundedness of $\lVert \mu_A \rVert$ and $\Psi(\Sigma_A)$ is a direct consequence of Corollary~\ref{cor:boundedness-sigma-mu} with $g = f$ and $B = A$, which implies by Lemma~\ref{lemma:equivalences} that $1/\lambda_1(\Sigma_A)$, $\lambda_d(\Sigma_A)$ and $\frob{\Sigma_A - I}$ are bounded. In turn, this implies the boundedness of $1/\lambda_1(\Sigma_\proj)$, $\lambda_d(\Sigma_\proj)$ and $\frob{\Sigma_\proj - I}$ by Lemma~\ref{lemma:proj-*} (applied with $\Sigma = \Sigma_A$, so that $\Sigma' = \Sigma_\proj$). Continuing, this implies the boundedness of $\Psi(\Sigma_\proj)$ by Lemma~\ref{lemma:equivalences}: thus, $\lVert \mu_A \rVert$ and $\Psi(\Sigma_\proj)$ are bounded, which implies by the same arguments that $g_\proj$ is efficient in high dimension for $A$ and that $1/\lambda_1(\Sigma_\proj)$, $\frob{\Sigma_\proj - I}$ and $\lambda_d(\Sigma_\proj)$ are bounded.
\end{proof}

It is clear that the arguments developed above apply when bounded is replaced with bounded whp. For the record, we state the generalization of the previous result that we will need later.

\begin{corollary} \label{cor:Sigma-proj-mu-A-bounded-stoch}
	For each $d$, let $B \subset \R^d$ be a random measurable set. If $1/p_f(B)$ is bounded whp, then $\lVert \mu_B \rVert$ and $\Psi(\Sigma_B)$ are bounded whp.
\end{corollary}

\subsection{High-dimensional efficiency of $\hat g_\proj$}

We first prove that $\lVert \hat \mu_A \rVert$ is bounded whp, and then that $\Psi(\hat \Sigma_\proj)$ is bounded whp.

\subsubsection{High-probability boundedness of $\lVert \hat \mu_A \rVert$} \label{subsub:tightness-hat-mu-A}

\begin{lemma} \label{lemma:tightness-hat-mu-A}
	We have
   	\begin{equation} \label{eq:bound-mu_A}
	    \E \left( \lVert \hat \mu_A - \mu_A \rVert^2 \right) \leq \frac{d}{\npar} \lambda_d(\Sigma_A).
    \end{equation}
	In particular, if $\inf_d p_f(A) > 0$ and $\npar \gg d$, then $\lVert \hat \mu_A - \mu_A \rVert \Rightarrow 0$ and $\lVert \hat \mu_A \rVert$ is bounded whp.
\end{lemma}

\begin{proof}
	Let us first prove~\eqref{eq:bound-mu_A}. Recall the definition~\eqref{eq:hat-mu-A-hat-Sigma-A} of $\hat \mu_A = \frac{1}{\npar} \sum_i Y_{A,i}$ with the $Y_{A,i}$'s i.i.d.\ distributed according to $f|_A$, so that
	\[ \E \left( \lVert \hat \mu_A - \mu_A \rVert^2 \right) = \frac{1}{\npar^2} \E \left( \sum_{i,j} (Y_{A,i} - \mu_A)^\top(Y_{A,j} - \mu_A) \right). \]
	Since the $Y_{A,i} - \mu_A$'s are i.i.d.\ and centered, we obtain
	\begin{align*}
		\E \left( \lVert \hat \mu_A - \mu_A \rVert^2 \right) & = \frac{1}{\npar} \E \left( (Y_{A,1} - \mu_A)^\top(Y_{A,1} - \mu_A) \right)\\
		& = \frac{1}{\npar} \E \left( \tr((Y_{A,1} - \mu_A) (Y_{A,1} - \mu_A)^\top) \right)
	\end{align*}
	which gives
	\[ \E \left( \lVert \hat \mu_A - \mu_A \rVert^2 \right) = \frac{1}{\npar} \tr \left( \Sigma_A \right) \]
	by commuting the trace and expectation operators. Since $\tr(\Sigma_A) \leq \lambda_d(\Sigma_A) d$ this gives~\eqref{eq:bound-mu_A}. Let us now assume that $\inf_d p_f(A)>0$ and $\npar \gg d$. Then $\lambda_d(\Sigma_A)$ is bounded by Corollary~\ref{cor:Sigma-proj-mu-A-bounded}, and so we obtain the result.
\end{proof}

\subsubsection{High-probability boundedness of $\Psi(\hat \Sigma_\proj)$} \label{subsub:tightness-hat-Sigma-proj}

To prove the fact that $\Psi(\hat \Sigma_\proj)$ is bounded whp, we need to study the spectrum of $\hat \Sigma_A$. Let in the sequel $\tilde Y_{A, i} = \Sigma_A^{-1/2} (Y_{A,i} - \mu_A)$ and $M$ be the $n \times d$ matrix with rows the $\tilde Y_{A,i}^\top$: then one can check that
\begin{equation} \label{eq:decomposition-S}
	\hat \Sigma_A = \Sigma_A^{1/2} \hat S \Sigma_A^{1/2} - (\hat \mu_A - \mu_A) (\hat \mu_A - \mu_A)^\top \ \text{ with } \ \hat S = \frac{1}{\npar} M^\top M.
\end{equation}
We will use results from~\cite{Vershynin10-0} in the area of non-asymptotic random matrix theory. The next lemma controls the subgaussian norm of the $\tilde Y_{A,i}$'s. According to the definitions~5.7 and 5.22 in~\cite{Vershynin10-0}, the subgaussian norm $\lVert Z \rVert_{\psi_2}$ of a $d$-dimensional random vector $Z$ is given by
\[ \lVert Z \rVert_{\psi_2} = \sup_{x: \lVert x \rVert = 1} \sup_{q \geq 1} q^{-1/2} \left( \E\lvert Z^\top x \rvert^q \right)^{1/q} = \sup_{x: \lVert x \rVert = 1} \lVert x^\top Z \rVert_{\psi_2}. \]
In the sequel, we denote by $Y_A$ and $\tilde Y_A$ random variables distributed as $Y_{A,i}$ and $\tilde Y_{A,i}$, respectively.

\begin{lemma}\label{lemma:subgaussian-norm}
	If $\inf_d p_f(A) > 0$, then $\sup_d \lVert \tilde Y_A \rVert_{\psi_2} < \infty$.
\end{lemma}

\begin{proof}
	Using the triangle inequality and the fact that the subgaussian norm of a constant vector is its norm, we obtain
 \[ \lVert \tilde Y_A \rVert_{\psi_2} = \lVert \Sigma_A^{-1/2} (Y_A - \mu_A) \rVert_{\psi_2} \leq \lVert \Sigma_A^{-1/2} Y_A \rVert_{\psi_2} + \lVert \Sigma_A^{-1/2} \mu_A \rVert. \]
	Note that $\lVert \Sigma_A^{-1/2} \mu_A \rVert = (\mu_A^\top \Sigma_A^{-1} \mu_A)^{1/2} \leq \lVert \mu_A \rVert / \lambda_1(\Sigma_A)^{1/2}$. Further, let $x \in \R^d$ with $\lVert x \rVert = 1$ and $Y \sim f$: then by definition of $Y_A$, for any $q \geq 1$ we have
 \[ \E\lvert x^\top \Sigma_A^{-1/2} Y_A \rvert^q = \E \left( \lvert x^\top \Sigma_A^{-1/2} Y \rvert^q\mid Y \in A \right) \leq \frac{1}{p_f(A)} \E \left( \lvert x^\top \Sigma_A^{-1/2} Y \rvert^q \right) \]
	and so (using $1/p_f(A)^{1/q} \leq 1/p_f(A)$ for $q \geq 1$)
	\[ \lVert x^\top \Sigma_A^{-1/2} Y_A \rVert_{\psi_2} \leq \frac{1}{p_f(A)} \lVert x^\top \Sigma_A^{-1/2} Y \rVert_{\psi_2}. \]
	For any centered Gaussian random variable $Z$, we have $\lVert Z \rVert_{\psi_2} \leq C \Var(Z)^{1/2}$ for some absolute constant $C$ (see for instance~\cite[Example 5.8]{Vershynin10-0}). Applying this to $Z = x^\top \Sigma^{-1/2}_A Y$, we obtain
	\[ \lVert x^\top \Sigma_A^{-1/2} Y \rVert_{\psi_2} \leq C \Var(x^\top \Sigma_A^{-1/2} Y)^{1/2} = C \sqrt{x^\top \Sigma_A^{-1} x} \leq \frac{C}{\lambda_1(\Sigma_A)^{1/2}}. \]
	Gathering the previous bounds, we therefore obtain
 \[ \lVert \Sigma_A^{-1/2} (Y_A - \mu_A) \rVert_{\psi_2} \leq \frac{1}{\lambda_1(\Sigma_A)^{1/2}} \left( \frac{C}{p_f(A)} + \lVert \mu_A \rVert \right). \]
	Under the assumption $\inf_d p_f(A) > 0$, Corollary~\ref{cor:Sigma-proj-mu-A-bounded} implies that this upper bound is bounded, which proves the result.
\end{proof}

\begin{lemma} \label{lemma:S}
	Let $\delta = \max(\lvert \lambda_1(\hat S) - 1 \rvert, \lvert \lambda_d(\hat S) - 1 \rvert)$. If $\inf_d p_f(A) > 0$ and $\npar \gg d$, then $\frac{\npar}{d} \delta^2$ is bounded whp. In particular, $\delta \Rightarrow 0$.
\end{lemma}

\begin{proof}
	By definition, the $\tilde Y_{A, i}$'s are i.i.d.\ centered random vectors. Moreover, they are isotropic, meaning that their covariance matrix is equal to the identity~\cite[Definition 5.19]{Vershynin10-0}, and they are subgaussian since their subgaussian norm is finite by Lemma~\ref{lemma:subgaussian-norm}. If $s_1$ and $s_d$ are the smallest and largest singular values of $M$, then Theorem 5.39 in~\cite{Vershynin10-0} implies that for any $t \geq 0$,
	\[ \P \left( \sqrt{\npar} - C' \sqrt{d} - t \leq s_1 \leq s_d \leq \sqrt{\npar} + C' \sqrt{d} + t \right) \geq 1 - 2 e^{-ct^2} \]
	where the constants $c$ and $C'$ only depend on the subgaussian norm of $\tilde Y_A$. But since $\sup_d \lVert \tilde Y_A \rVert_{\psi_2} < \infty$ by Lemma~\ref{lemma:subgaussian-norm}, it follows that the constants $c$ and $C'$ can be chosen independent of $d$. Moreover, since $\hat S = \frac{1}{\npar} M^\top M$, we have
	\[ \lambda_1(\hat S) = \frac{1}{\npar} s_1^2 \ \text{ and } \ \lambda_d(\hat S) = \frac{1}{\npar} s_d^2, \]
	and so for $t = \sqrt d$, we obtain
	\begin{multline*}
		\P \left( \left( 1 - (C'+1) \sqrt{d/\npar} \right)^2 \leq \lambda_1(\hat S) \leq \lambda_d(\hat S) \leq \left( 1 + (C'+1) \sqrt{d/\npar} \right)^2 \right)\\
		\geq 1 - 2 e^{-c d}
	\end{multline*}
	From there, one can easily derive the result through elementary manipulation.
\end{proof}

\begin{corollary} \label{cor:tightness}
	If $\inf_d p_f(A) > 0$ and $\npar \gg d$, then the sequences $1/\lambda_1(\hat \Sigma_\proj)$ and $\lambda_d(\hat \Sigma_\proj)$ are bounded whp.
\end{corollary}

\begin{proof}
	Lemma~\ref{lemma:proj-*} with $\Sigma = \hat \Sigma_A$ and $d_k = d_k$ (so that $\Sigma' = \hat \Sigma_\proj$) give
	\[ \lambda_1(\hat \Sigma_\proj) \geq \min(1, \lambda_1(\hat \Sigma_A)), \ \lambda_d(\hat \Sigma_\proj) \leq \max(1, \lambda_d(\hat \Sigma_A)) \]
	and $\frob{\hat \Sigma_\proj - I} \leq \frob{\hat \Sigma_A - I}$. Therefore, it is enough to show that $1/\lambda_1(\hat \Sigma_A)$ and $\lambda_d(\hat \Sigma_A)$ are bounded whp. Since $\delta \Rightarrow 0$ by Lemma~\ref{lemma:S}, we have $\lambda_1(\hat S) \Rightarrow 1$ and $\lambda_d(\hat S) \Rightarrow 1$, and so the sequences $1/\lambda_1(\hat S)$ and $\lambda_d(\hat S)$ are bounded whp. Thus, we only have to transfer this result to $\hat \Sigma_A$. Let $x \in \R^d$ with $\lVert x \rVert = 1$, and $y = \Sigma_A^{1/2} x$: then by definition (see~\eqref{eq:decomposition-S}), we have
	\[ x^\top \hat \Sigma_A x = y^\top \hat S y - (x^\top (\hat \mu_A - \mu_A))^2. \]
	In particular,
	\[ x^\top \hat \Sigma_A x \leq \lambda_d(\hat S) \lVert y \rVert^2 = \lambda_d(\hat S) x^\top \Sigma_A x \leq \lambda_d(\hat S) \lambda_d(\Sigma_A) \]
	and so
	\[ \lambda_d(\hat \Sigma_A) \leq \lambda_d(\hat S) \lambda_d(\Sigma_A). \]
	Since $\lambda_d(\hat S)$ is bounded whp and $\lambda_d(\Sigma_A)$ is bounded, this that $\lambda_d(\hat \Sigma_A)$ is bounded whp. We show that $1/\lambda_1(\hat \Sigma_A)$ is bounded whp with similar arguments: we have
	\[ x^\top \hat \Sigma_A x \geq \lambda_1(\hat S) x^\top \Sigma_A x - \lVert \hat \mu_A - \mu_A \rVert^2 \geq \lambda_1(\hat S) \lambda_1(\Sigma_A) - \lVert \hat \mu_A - \mu_A \rVert^2 \]
	and so
	\[ \lambda_1(\hat \Sigma_A) \geq \lambda_1(\hat S) \lambda_1(\Sigma_A) - \lVert \hat \mu_A - \mu_A \rVert^2. \]
	Since $\lVert \hat \mu_A - \mu_A \rVert \Rightarrow 0$ when $n_g\gg d$ by Lemma~\ref{lemma:tightness-hat-mu-A}, $1/\lambda_1(\hat S)$ is bounded whp by Lemma~\ref{lemma:S} and $1/\lambda_1(\Sigma_A)$ is bounded by Corollary~\ref{cor:Sigma-proj-mu-A-bounded}, the previous inequality gives that $1/\lambda_1(\hat \Sigma_A)$ is bounded whp.
\end{proof}

\begin{lemma}
	If $\inf_d p_f(A) > 0$ and $\npar \gg rd$, then $\Psi(\hat \Sigma_\proj)$ is bounded whp.
\end{lemma}

\begin{proof}
	According to Corollary~\ref{cor:tightness}, $1/\lambda_1(\hat \Sigma_\proj)$ and $\lambda_d(\hat \Sigma_\proj)$ are bounded whp. Thus, in order to show that $\Psi(\hat \Sigma_\proj)$ is bounded whp, it remains to show in view of Lemma~\ref{lemma:equivalences} that $\frob{\hat \Sigma_\proj - I}$ is bounded whp. Define
	\[ \Sigma'_\proj = \sum_{k=1}^r (v_k - 1) d_k d_k^\top + I \ \text{ with } \ v_k = d_k^\top \Sigma_A d_k. \]
 According to Lemma~\ref{lemma:proj-*}, we have that $\frob{\Sigma'_\proj - I} \leq \frob{\Sigma_A - I}$. Since $\frob{\Sigma_A - I}$ is bounded by Corollary~\ref{cor:Sigma-proj-mu-A-bounded}, we obtain that $\frob{\Sigma'_\proj - I}$ is bounded. By the triangle inequality, it is therefore enough to prove that $\frob{\hat \Sigma_\proj - \Sigma'_\proj}$ is bounded whp. By definition we have
	\[ \hat \Sigma_\proj - \Sigma'_\proj = \sum_{k=1}^r (\hat v_k - v_k) d_k d_k^\top. \]
	Therefore, the eigenvalues of $\hat \Sigma_\proj - \Sigma'_\proj$ are the $\hat v_k - v_k$ for $k = 1, \ldots, r$, and $0$ with multiplicity $d-r$. Since $\hat \Sigma_\proj - \Sigma'_\proj$ is symmetric, the square of its Frobenius norm is equal to the sum of the square of its eigenvalues, and since at most $r$ of them are non-zero, we obtain
	\[ \frob{ \hat \Sigma_\proj - \Sigma'_\proj}^2 \leq r \varepsilon \ \text{ with } \ \varepsilon = \max \left( \lambda_1(\hat \Sigma_A - \Sigma_A)^2, \lambda_d(\hat \Sigma_A - \Sigma_A)^2 \right). \]
	By~\eqref{eq:decomposition-S} we have
	\[ \hat \Sigma_A - \Sigma_A = \Sigma_A^{1/2} (\hat S - I) \Sigma_A^{1/2} - (\hat \mu_A - \mu_A) (\hat \mu_A - \mu_A)^\top \]
	and so if we let $\delta = \max(\lvert \lambda_1(\hat S) - 1 \rvert, \lvert \lambda_d(\hat S) - 1 \rvert)$ as in Lemma~\ref{lemma:S}, we obtain for any $x \in \R^d$ with $\lVert x \rVert = 1$
	\[ \left \lvert x^\top (\hat \Sigma_A - \Sigma_A) x \right \rvert \leq \lambda_d(\Sigma_A) \delta + \lVert \hat \mu_A - \mu_A \rVert^2. \]
	By definition of $\varepsilon$ and the variational characterization of eigenvalues, this implies that $\varepsilon^{1/2} \leq \lambda_d(\Sigma_A) \delta + \lVert \hat \mu_A - \mu_A \rVert^2$ and since $\frob{ \hat \Sigma_\proj - \Sigma'_\proj}^2 \leq r \varepsilon$, we finally get
	\[ \frob{ \hat \Sigma_\proj - \Sigma'_\proj}^2 \leq r \left( \lambda_d(\Sigma_A) \delta + \lVert \hat \mu_A - \mu_A \rVert^2 \right)^2. \]
	Given that $\lambda_d(\Sigma_A)$ is bounded by Corollary~\ref{cor:Sigma-proj-mu-A-bounded}, the proof will be complete if we prove that $r \delta^2 \Rightarrow 0$ and $r \lVert \hat \mu_A - \mu_A \rVert^4 \Rightarrow 0$, which is what we do in the rest of the proof.
	
	The fact that $r \delta^2 \Rightarrow 0$ is a direct consequence of Lemma~\ref{lemma:S}, which implies that $\P(r \delta^2 \leq C (rd/\npar)) \to 1$ (which gives $r \delta^2 \Rightarrow 0$ since $r d / \npar \to 0$). On the other hand,~\eqref{eq:bound-mu_A} directly implies that $r \lVert \hat \mu_A - \mu_A \rVert^2 \Rightarrow 0$ when $\npar \gg rd$, which implies that $r \lVert \hat \mu_A - \mu_A \rVert^4 \Rightarrow 0$. The proof is therefore complete.
\end{proof}

\section{Proof of Theorem~\ref{thm:p}} \label{sec:proof-p}

From~\eqref{eq:D-g-mid-B-g'-tmp}, one can derive the following two identities:
\begin{equation} \label{eq:D-g*-g-opt}
	D(f|_A || g_A) = -\log p_f(A) - \Psi(\Sigma_A) - \frac{1}{2} \lVert \mu_A \rVert^2
\end{equation}
and
\begin{equation} \label{eq:D-f-g-opt}
	D(f || g_A) = \Psi(\Sigma_A^{-1}) + \frac{1}{2} \mu_A^\top \Sigma_A^{-1} \mu_A.
\end{equation}
Assume now that $p_f(A) \to 0$ and that $\sup_d D(f || g_A) < \infty$: in order to prove Theorem~\ref{thm:p}, it is enough to prove that $D(f|_A || g_A) \to \infty$. In view of~\eqref{eq:D-f-g-opt}, the boundedness of $D(f || g_A)$ implies that of $\Psi(\Sigma_A^{-1})$ and of $\mu_A^\top \Sigma_A^{-1} \mu_A$. The boundedness of $\Psi(\Sigma_A^{-1})$ implies by Lemma~\ref{lemma:equivalences} that of $\Psi(\Sigma_A)$ and of $\lambda_d(\Sigma_A)$. Since
\[ \mu_A^\top \Sigma_A^{-1} \mu_A \geq \frac{\lVert \mu_A \rVert^2}{\lambda_d(\Sigma_A)}, \]
this implies the boundedness of $\lVert \mu_A \rVert$. Thus, we have proved that the sequences $\Psi(\Sigma_A)$ and $\lVert \mu_A \rVert$ are bounded: but $-\log p_f(A) \to +\infty$, and so $D(f|_A || g_A) \to \infty$ in view of~\eqref{eq:D-g*-g-opt} which proves the result.

\section{Proof of Theorem~\ref{thm:CE}} \label{sec:proof-CE}

\subsection{High-dimensional efficiency of $g_t$} \label{sub:g_t}

Let us now turn to the high-dimensional efficiency of $g_t$. We use throughout the notation introduced before Theorem~\ref{thm:CE}. We proceed by induction on $t \geq 0$, working with the following induction hypothesis.

\begin{HI-det}
	For $t \geq 0$, $\Psi(\Sigma_t)$, $\lVert \mu_t \rVert$ and $1/p_f(A_t)$ are bounded.
\end{HI-det}

Note that if $\Psi(\Sigma_t)$ and $\lVert \mu_t \rVert$ are bounded, then $g_t = N(\mu_t, \Sigma_t)$ is efficient in high dimension by Proposition~\ref{prop:first-step}. The additional requirement that $1/p_f(A_t)$ is bounded is here to pass through the induction.

\begin{lemma} \label{lemma:initialization}
	If for every $d$, $\varphi$ has no atom and $\inf_d \rho > 0$, then the deterministic induction hypothesis holds for $t = 0$.
\end{lemma}

\begin{proof}
    Let $F(x) = \P_f(\varphi(X) \leq x)$. Since $g_0 = f$ we have by definition of $A_0 = \{x: \varphi(x) > q_0\}$ and $q_0 = F^{-1}(1-\rho)$
    \[ p_f(A_0) = \P_{g_0}(X \in A_0) = \P_{g_0}(\varphi(X) > q_0) = 1 - F(F^{-1}(1-\rho)) = \rho \]
    using Lemma~\ref{lemma:no-atom} for the last equality. Since $\Sigma_0 = I$ and $\mu_0 = 0$ and we assume $\inf_d \rho > 0$, we get that $\Psi(\Sigma_0)$, $\lVert \mu_0 \rVert$ and $1/p_f(A_0)$ are bounded, i.e., the deterministic induction hypothesis holds for $t = 0$.
\end{proof}

We now prove the induction.

\begin{lemma} \label{lemma:induction-det}
	Assume that for every $d$, $\varphi$ has no atom and that $\inf_d \rho > 0$. If the deterministic induction hypothesis holds for some $t \geq 0$, then it holds at $t+1$.
\end{lemma}

\begin{proof}
	Assume that the deterministic induction hypothesis holds for some $t \geq 0$, i.e., $\Psi(\Sigma_t)$, $\lVert \mu_t \rVert$ and $1/p_f(A_t)$ are bounded, and let us show that this continues to hold for $t+1$. The boundedness of $1/p_f(A_t)$ implies by Corollary~\ref{cor:boundedness-sigma-mu} with $g = f$ and $B = A_t$ that $\Psi(\Sigma_{A_t})$ and $\lVert \mu_{A_t} \rVert$ are bounded. Since $\mu_{t+1} = \mu_{A_t}$ and $\Sigma_{t+1} = \Sigma_{A_t}$, it remains to prove that $1/p_f(A_{t+1})$ is bounded. Using Corollary~\ref{cor:jensen} with $B = A_{t+1}$ and $g = g_{t+1}$, we obtain
	\begin{equation} \label{eq:p-f-g}
		p_f(A_{t+1}) \geq p_{g_{t+1}}(A_{t+1}) \exp \left( - \Psi(\Sigma^{g_{t+1}}_{A_{t+1}}) - \frac{1}{2} \lVert \mu^{g_{t+1}}_{A_{t+1}} \rVert^2 \right).
	\end{equation}
	Recall that by definition of the CE scheme and Lemma~\ref{lemma:no-atom}, we have
	\[ p_{g_{t+1}}(A_{t+1}) = \P_{g_{t+1}}(\varphi(X) > q_{t+1}) = 1 - F(F^{-1}(1-\rho)) = \rho. \]
	Since we assume $\inf_d \rho > 0$, it remains only in view of~\eqref{eq:p-f-g} to prove that $\Psi(\Sigma^{g_{t+1}}_{A_{t+1}})$ and $\lVert \mu^{g_{t+1}}_{A_{t+1}} \rVert$ are bounded. But since $\lVert \mu_{t+1} \rVert$, $\Psi(\Sigma_{t+1})$ and $1/p_{g_{t+1}}(A_{t+1})$ are bounded, this follows precisely from Corollary~\ref{cor:boundedness-sigma-mu} with $g = g_{t+1}$ and $B = A_{t+1}$. Thus, the deterministic induction hypothesis holds at $t+1$.
\end{proof}

We can now prove the part of Theorem~\ref{thm:CE} that relates to $g_t$.

\begin{prop} \label{prop:HDE-g_t}
	If $1/p_f(A)$ and $1/\rho$ are bounded, and if for every $d$, $\varphi$ has no atom, then for every $t \geq 0$, $g_t$ is efficient in high dimension for $A$.
\end{prop}

\begin{proof}
	Combining Lemmas~\ref{lemma:initialization} and~\ref{lemma:induction-det}, we get that $\lVert \mu_t \rVert$ and $\Psi(\Sigma_t)$ are bounded for every $t \geq 0$. Combined with the assumption $\inf_d p_f(A) > 0$, this gives the result in view of Proposition~\ref{prop:first-step}.
\end{proof}

\subsection{High-dimensional efficiency of $\hat g_t$}

\subsubsection{Proof outline} \label{subsub:proof-outline}

Compared to $\hat g_\proj$, analyzing the CE scheme (i.e., showing that $\hat g_t$ is efficient in high-dimension) entails one significant additional difficulty which imposes the implicit growth rate $n \gg d^\kappa$ in Theorem~\ref{thm:CE}. In order to illustrate this difficulty, consider
\[ \hat \mu'_{t+1} = \frac{1}{\npar p_f(\hat A_t)} \sum_{i=1}^\npar \ell(Y_i) \xi_{\hat A_t}(Y_i) Y_i \ \text{ with } \ \ell = f/\hat g_t. \]
Compared to $\hat \mu_{t+1}$ in~\eqref{eq:CE-p-mu}, we have just replaced $\hat p_t$ by $p_f(\hat A_t)$, but thanks to this mild modification, we can use the CD bound~\eqref{eq:bound-CD}, conditional on $\hat g_t$ and $\hat A_t$, on every coordinate $k=1, \ldots, d$ with $\phi(x) = x(k) \xi_{\hat A_t}(x) / p_f(\hat A_t )$, to get a bound on $\lvert \hat \mu'_{t+1} - \mu_{\hat A_t} \rvert$. We will see below that this approach leads to a bound of the form
\[ \widehat \E \left( \lone{\hat \mu'_{t+1} - \mu_{\hat A_t}} \right) \leq \frac{Z'}{p_f(\hat A_t)} d n^{-\alpha/4} \]
with $Z'$ bounded whp, see Lemma~\ref{lemma:E-Z} and Lemma~\ref{lemma:CD-bounds} below ($\widehat \E$ will be introduced below also). What is important is that this bound holds for any $\alpha < \alpha_*(\hat \Sigma_t)$ (recall the definition~\eqref{eq:alpha*} of $\alpha_*$).

Thus, if we want to make this bound vanish (which is the first step toward the control of $\mu_{t+1}$), we need $d n^{-\alpha/4} \to 0$ for some $\alpha < \alpha_*(\hat \Sigma_t)$, i.e., $n \gg d^\kappa$ for some $\kappa > 4/\alpha_*(\hat \Sigma_t)$. This approach ultimately gives a control on $\hat \mu_{t+1}$, but at the expense of a \textit{random} growth rate for $n$, which is unsatisfactory. As discussed in the end of Section~\ref{sub:CE}, the intuition $\hat \Sigma_t \approx \Sigma_t$ suggests to try and show that $\alpha_*(\hat \Sigma_t) \approx \alpha_*(\Sigma_t)$, which is tantamount to showing that $\lambda_1(\hat \Sigma_t) \approx \lambda_1(\Sigma_t)$. However, controlling smallest eigenvalues of random matrices is a difficult problem, and it seems that justifying the approximation $\lambda_1(\hat \Sigma_t) \approx \lambda_1(\Sigma_t)$ would require additional technical assumptions, e.g., on the growth rate of $m$ and regularly properties for $\varphi$. Here we adopt a different approach, and just prove the existence of $\underline \alpha > 0$ such that $\P(\alpha_*(\hat \Sigma_t) \geq \underline \alpha) \to 1$. The approach outlined above then provides a control of $\hat \mu_{t+1}$ provided $\npar \gg d^{4/\underline \alpha}$.

As in the control of $g_t$, the control of $\hat g_t$ proceeds by induction. To that purpose we need the following stochastic version of the previous deterministic induction hypothesis.

\begin{HI-stoch}
	Let $t \geq 0$. We say that the stochastic induction hypothesis holds at time $t$ if $\Psi(\hat \Sigma_t)$, $\lVert \hat \mu_t \rVert$ and $1/p_f(\hat A_t)$ are bounded whp.
\end{HI-stoch}

The initialization of the induction will be carried out in Lemma~\ref{lemma:init-induction-hat-g_t}, while the induction itself is treated in Theorem~\ref{thm:CE-induction}.

\begin{lemma} \label{lemma:init-induction-hat-g_t}
	Assume that:
	\begin{itemize}
		\item $\inf_d \rho > 0$;
		\item for every $d$, $\varphi$ has no atom;
		\item $m \to \infty$.
	\end{itemize}
	Then for any $t \geq 0$, $1/p_{\hat g_t}(\hat A_t)$ is bounded whp. In particular, the stochastic induction hypothesis holds at time $t = 0$.
\end{lemma}

\begin{thm}\label{thm:CE-induction}
	Assume that:
	\begin{itemize}
		\item $\inf_d p_f(A) > 0$;
		\item $\inf_d \rho > 0$;
		\item for every $d$, $\varphi$ has no atom;
		\item $m \to \infty$.
	\end{itemize}
     Under these assumptions, if the stochastic induction hypothesis holds at some time $t \geq 0$, then there exists a constant $\kappa > 0$ such that if $n \gg d^\kappa$, then the stochastic induction hypothesis holds at time $t+1$.
\end{thm}

Before proceeding to the proof of this result, let us complete the proof of the part of Theorem~\ref{thm:CE} related to $\hat g_t$ based on Lemma~\ref{lemma:init-induction-hat-g_t} and Theorem~\ref{thm:CE-induction}.

\begin{prop}
	Assume that:
	\begin{itemize}
		\item $\inf_d p_f(A) > 0$;
		\item $\inf_d \rho > 0$;
		\item for every $d$, $\varphi$ has no atom;
		\item $m \to \infty$.
	\end{itemize}
     Then for any $t \geq 0$ there exists a constant $\kappa_t > 0$ such that if $\npar \gg d^{\kappa_t}$, then $\hat g_t$ is efficient in high dimension for $A$.
\end{prop}

\begin{proof} [Proof based on Lemma~\ref{lemma:init-induction-hat-g_t} and Theorem~\ref{thm:CE-induction}]
	Lemma~\ref{lemma:init-induction-hat-g_t} implies that the stochastic induction hypothesis holds at time $0$, and Theorem~\ref{thm:CE-induction} then implies that it holds for every $t \geq 0$. Thus, $\lVert \hat \mu_t \rVert$ and $\Psi(\hat \Sigma_t)$ are bounded whp, and $1/p_f(A)$ is bounded: Proposition~\ref{prop:first-step} then implies that $\hat g_t$ is efficient in high dimension for~$A$.
\end{proof}

\subsubsection{Control of $p_{\hat g_t}(\hat A_t)$ and induction initialization}

In this section we prove Lemma~\ref{lemma:init-induction-hat-g_t}.

\begin{lemma} \label{lemma:quantile}
	For each $d$, let:
	\begin{itemize}
		\item $U_1, \dots, U_m$ be $m = m(d)$ i.i.d.\ real-valued random variables with cumulative distribution function $F(u) = \P(U \leq u)$, that may depend on $d$;
		\item $F$ be continuous;
		\item $\varrho \in (0,1)$ and $q = F^{-1}(1-\varrho)$;
		\item $\hat q = U_{([(1-\varrho) m])}$ the empirical estimation of $q$, with $U_{(1)} \leq \cdots \leq U_{(m)}$.
	\end{itemize}
	 Assume that $m \to \infty$ and $\inf_d \varrho > 0$. Then $(1-F(\hat q)) / \varrho \Rightarrow 1$, and in particular, $1/(1-F(\hat q))$ is bounded whp.
\end{lemma}

\begin{proof}
	We have
	\[ \P((1-F(\hat q)) / \varrho \leq x) = \P(F(\hat q) \geq 1-\varrho x) = \P(\hat q \geq F^{-1}(1-\varrho x))  \]
	with the second equality coming from the fact that $F^{-1}$ is the left-continuous inverse, so $F(x) \geq t \Leftrightarrow x \geq F^{-1}(t)$. Let $n = [(1-\varrho)m]$: by definition of $\hat q$, we have
	\[ \P(\hat q \geq F^{-1}(1-\varrho x)) = \P(U_{(n)} \geq F^{-1}(1-\varrho x)). \]
	Since $U_{(k)}$ is the $k$th largest sample among the $U_i$'s, we have
	\[ U_{(n)} \geq F^{-1}(1-\varrho x) \Longleftrightarrow \#\{i: U_i \geq F^{-1}(1-\varrho x)\} \geq m-n+1. \]
	Since the $U_i$'s are i.i.d., the random variable $\#\{i: U_i \geq F^{-1}(1-\varrho x)\}$ follows a binomial distribution with parameters $m$ and
	\[ \P(U_1 \geq F^{-1}(1-\varrho x)) = 1 - F(F^{-1}(1-\varrho x)) = \varrho x, \]
	both equalities coming from the fact that $F$ is continuous. Thus if $B_{q}$ denotes a binomial random variable with parameters $m$ and $q$, we obtain
	\[ \P(\hat q \geq F^{-1}(1-\varrho x)) = \P(B_{\varrho x} \geq m - n + 1). \]
	By considering Laplace transforms, one easily sees that $B_{\varrho x} / (\varrho m) \Rightarrow x$. Since $(m-n+1)/(\varrho m) \to 1$ (using $\inf_d \varrho > 0$), we obtain
	\[ \P(\hat q \geq F^{-1}(1-\varrho x)) \to \Indicator{x \geq 1} \]
	for $x \neq 1$, which implies the desired convergence $(1-F(\hat q)) / \varrho \Rightarrow 1$. As this clearly implies $1/(1-F(\hat q))$ that is bounded whp, this concludes the proof.
\end{proof}

\begin{proof} [Proof of Lemma~\ref{lemma:init-induction-hat-g_t}]
Let $F(u) = \P_{\hat g_t}(\varphi(X) \leq u)$: by definition of $\hat A_t$ we have $p_{\hat g_t}(\hat A_t) = 1-F(\hat q_t)$. Moreover, since $\varphi$ has no atom, Lemma~\ref{lemma:no-atom} implies that~$F$ is continuous, and so Lemma~\ref{lemma:quantile} (applied with $U_k = \varphi(Y'_k)$) implies that $p_{\hat g_t}(\hat A_t) / \rho \Rightarrow 1$, where the convergence holds under $\P_{\hat g_t}$, so conditionally on $\hat g_t$. In particular, for $x > 1/\inf_d \rho$ we have
	\[ \P_{\hat g_t}(1/p_{\hat g_t}(\hat A_t) \geq x) = \P_{\hat g_t}(p_{\hat g_t}(\hat A_t)/\rho \leq 1/(\rho x)) \Rightarrow 0 \]
	and so $\P(1/p_{\hat g_t}(\hat A_t) \geq x) \to 0$ as well, by the bounded convergence theorem, which proves that $p_{\hat g_t}(\hat A_t)$ is bounded whp.
 
 Concerning the stochastic induction hypothesis at time $t = 0$, note that for $t = 0$ we have $\hat \mu_0 = \mu_0 = 0$ and $\hat \Sigma_0 = I = \Sigma_0$, which readily entails that $\Psi(\hat \Sigma_0)$ and $\lVert \hat \mu_0 \rVert$ are bounded. Further, since $\hat g_0 = f$ we have $1/p_f(\hat A_0) = 1/p_{\hat g_0}(\hat A_0)$ which was just proved to be bounded whp.
\end{proof}

\subsubsection{Additional notation and preliminary results}

Before proceeding to the proof of Theorem~\ref{thm:CE-induction}, let us establish some preliminary results and introduce additional notation.

\begin{lemma} \label{lemma:alpha_*}
	Assume that the stochastic induction hypothesis holds at time $t$. Then $D(f || \hat g_t)$ and $1/\alpha_*(\hat \Sigma_t)$ are bounded whp. In particular, there exists $\underline \alpha > 0$ such that the event $\calE$ defined by
\[ \calE = \{D(f || \hat g_t) \leq \log n\} \cap \{ \underline \alpha < \alpha_*(\hat \Sigma_t) \} \]
holds with high probability, i.e., $\P(\calE) \to 1$.\\
\end{lemma}

\begin{proof}
	From~\eqref{eq:D-f-g}, we obtain
	\[ D(f || \hat g_t) = \Psi(\hat \Sigma_t^{-1}) + \frac{1}{2} \hat \mu_t^\top \hat \Sigma_t^{-1} \hat \mu_t \leq \Psi(\hat \Sigma_t^{-1}) + \frac{1}{2 \lambda_1(\hat \Sigma_t)} \lVert \hat \mu_t \rVert^2. \]
	Since $\Psi(\hat \Sigma_t)$ and $\lVert \hat \mu_t \rVert$ are bounded whp by assumption, Lemma~\ref{lemma:equivalences} implies that $D(f || \hat g_t)$ is bounded whp by the inequality of the previous display. Moreover, since $\Psi(\hat \Sigma_t)$ is bounded whp by the stochastic induction hypothesis, this implies that $1/\lambda_1(\hat \Sigma_t)$ is bounded whp by Lemma~\ref{lemma:equivalences}, which implies that $1/\alpha_*(\hat \Sigma_t)$ is bounded whp by definition of $\alpha_*$ in \eqref{eq:alpha*}.
\end{proof}

In the sequel, we assume that the stochastic induction hypothesis holds at time $t$. We fix a constant $\underline \alpha$ given by the previous lemma and we consider the event $\calE$ defined there. Let in the sequel
\[ \widehat \P = \P(\, \cdot \mid \hat g_t, \hat A_t, \calE) \]
be the random distribution conditional on $\hat g_t$, $\hat A_t$ and the event $\calE$. The motivation for introducing $\widehat \P$ is that conditioning $\hat g_t$, $\hat A_t$ and the event $\calE$ will allow us to use the CD bound~\eqref{eq:bound-CD}.

We consider an additional constant $\alpha < \underline \alpha$, and we define $Z = 0$ if $\alpha_*(\hat \Sigma_t) \leq \underline \alpha$, and
\begin{equation} \label{eq:Z}
	Z = \exp \left( \alpha D(f || \hat g_t) + \frac{1}{2} q \alpha^2 \lVert \hat \Sigma_t^{-1} \hat \mu_t \rVert^2 + \frac{\alpha}{2\underline\alpha} \Psi((\underline\alpha+1)I - \underline\alpha \hat \Sigma_t^{-1}) \right)
\end{equation}
if $\underline \alpha < \alpha_*(\hat \Sigma_t)$, with $q = \underline\alpha / (\underline\alpha - \alpha)$. Note that $Z$ is the bound~\eqref{eq:tail-exp} of the $\alpha$th-moment of the likelihood ratio between $f$ and $\hat g_t$. We also define
\[ Z' = 3 e^{\alpha D(f || \hat g_t)} Z^{1/2}. \]
We will use the following result on $Z$ and $Z'$.

\begin{lemma} \label{lemma:E-Z}
	If the stochastic induction hypothesis holds at time $t$, then $Z$ and $Z'$ are bounded whp.
\end{lemma}

\begin{proof}
	Recall that $Z' = 3 e^{\alpha D(f || \hat g_t)} Z^{1/2}$ with $Z$ defined in~\eqref{eq:Z}: since $D(f || \hat g_t)$, $\Psi(\hat \Sigma_t)$ and $\lVert \hat \mu_t \rVert$ are bounded whp by the stochastic induction hypothesis and Lemma~\ref{lemma:alpha_*}, it is enough in view of~\eqref{eq:Z} to show that $\Psi((\alpha+1) I - \alpha \hat \Sigma_t^{-1})$ is bounded whp. For $i \in \{1,d\}$, we have
	\begin{align*}
		\lambda_i \left( (\alpha+1) I - \alpha \hat \Sigma_t^{-1} \right) & = \alpha + 1 + \lambda_i(- \alpha \hat \Sigma^{-1}_t)\\
		& = \alpha + 1 - \alpha \lambda_{d+1-i}(\hat \Sigma^{-1}_t)\\
		& = \alpha + 1 - \frac{\alpha}{\lambda_i(\hat \Sigma_t)}\\
		& = 1 - \frac{1 - \lambda_i(\hat \Sigma_t)}{\lambda_i(\hat \Sigma_t)} \alpha. 
	\end{align*}
	Since $\Psi(\hat \Sigma_t)$ is bounded whp by the stochastic induction hypothesis, $\lambda_d(\hat \Sigma_t)$ and $1/\lambda_1(\hat \Sigma_t)$ are bounded whp by Lemma~\ref{lemma:equivalences}, and so the previous display implies that $1/\lambda_1((\alpha+1) I - \alpha \hat \Sigma_t^{-1})$ and $\lambda_d((\alpha+1) I - \alpha \hat \Sigma_t^{-1})$ are also bounded whp. Moreover,
	\[ \frob{(\alpha+1) I - \alpha \hat \Sigma_t^{-1} - I} = \alpha \frob{\hat \Sigma^{-1}_t - I} \]
	which is bounded whp, again as a consequence of the assumption that $\Psi(\hat \Sigma_t)$ is bounded whp and Lemma~\ref{lemma:equivalences}. Invoking Lemma~\ref{lemma:equivalences}, we obtain $\Psi((\alpha+1) I - \alpha \hat \Sigma_t^{-1})$ which concludes the proof.
\end{proof}

\subsubsection{Induction} \label{subsub:tightness}

We now prove Theorem~\ref{thm:CE-induction} that the induction goes through. So in the rest of this section, we assume that the assumptions of Theorem~\ref{thm:CE-induction} hold: in particular, the stochastic induction hypothesis holds at time $t$. We identify growth rates for $n$ that guarantee that $\Psi(\hat \Sigma_{t+1})$, $\lVert \hat \mu_{t+1} \rVert$ and $1/p_f(\hat A_{t+1})$ are bounded whp. We begin with the following lemma, which follows by combining the CD bound~\eqref{eq:bound-CD} and the bound~\eqref{eq:tail-exp} on the exponential moments of the log-likelihood. In the sequel, define
\[ \ell = \frac{f}{\hat g_t}. \]

\begin{lemma}\label{lem:}
	Let $Y_i$ i.i.d.~$\sim \hat g_t$, $Y \sim f$, $d' \in \N\setminus\{0\}$, and $\phi: \R^d \to \R^{d'}$ measurable written $\phi(x) = (\phi_1(x), \ldots, \phi_{d'}(x))$ for $x \in \R^d$, with $\phi_k: \R^d \to \R$. Then
	\begin{equation} \label{eq:CD+exp-tail}
		\widehat \E \left( \Lone{\frac{1}{\npar} \sum_{i=1}^{\npar} \ell(Y_i) \phi(Y_i) - \widehat \E(\phi(Y))} \right) \leq Z' \left ( \sum_{k=1}^{d'} \sqrt{\widehat \E\left( \phi_k(Y)^2 \right)} \right) n^{-\alpha/4}.
	\end{equation}
\end{lemma}

\begin{proof}
	Since $\widehat \P(\calE) = 1$, we can use~\eqref{eq:bound-CD} with $g = \hat g_t$ and $\phi = \phi_k$ to obtain
	\begin{multline*}
		\widehat \E \left( \left \lvert \frac{1}{n} \sum_{i=1}^n \ell(Y_i) \phi_k(Y_i) - \widehat \E(\phi_k(Y)) \right \rvert \right) \leq \left( \widehat \E \left(\phi_k(Y)^2 \right) \right)^{1/2} \times\\
		\left[ \left(\frac{e^{D(f || \hat g_t)}}{n}\right)^{1/4} + 2 \left( \widehat \P \left(L(Y) \geq \frac{1}{2} \log n + \frac{1}{2} D(f || \hat g_t) \right) \right)^{1/2} \right]
	\end{multline*}
	with $L = \log(f/\hat g_t)$. Concerning the tail of the log-likelihood, we have
	\begin{align*}
		\widehat \P \left(L(Y) \geq \frac{1}{2} \log n + \frac{1}{2} D(f || \hat g_t) \right) & = \widehat \P \left( e^{\alpha L(Y)} \geq (ne^{D(f || \hat g_t)})^{\alpha/2} \right)\\
		& \leq (ne^{D(f || \hat g_t)})^{-\alpha/2} \widehat \E \left( e^{\alpha L(Y)} \right)\\
		& \leq (ne^{D(f || \hat g_t)})^{-\alpha/2} Z
	\end{align*}
	using Lemma~\ref{lemma:tail-exp} for the last inequality (which we can invoke, since by definition of $\alpha, \underline \alpha$ and $\widehat \P$ we have $\widehat \P(\alpha < \underline\alpha < \alpha_*(\hat \Sigma_t)) = 1$). This leads to
	\begin{multline*}
		\widehat \E \left( \left \lvert \frac{1}{n} \sum_{i=1}^n \ell(Y_i) \phi_k(Y_i) - \widehat \E(\phi_k(Y)) \right \rvert \right) \leq \left( \widehat \E \left(\phi_k(Y)^2 \right) \right)^{1/2} \times\\
		\left[ \left(\frac{e^{D(f || \hat g_t)}}{n}\right)^{1/4} + 2 \left( \frac{e^{D(f || \hat g_t)}}{n} \right)^{\alpha/4} Z^{1/2} \right].
	\end{multline*}
	Since $e^{D(f||\hat g_t)}/n \leq 1$ (since we are in the event $\{D(f || \hat g_t) \leq \log n\}$) and $\alpha < 1$, we have
	\[ \left(\frac{e^{D(f || \hat g_t)}}{n}\right)^{1/4} \leq \left( \frac{e^{D(f || \hat g_t)}}{n} \right)^{\alpha/4} \]
	and so we get
	\begin{multline*}
		\widehat \E \left( \left \lvert \frac{1}{n} \sum_{i=1}^n \ell(Y_i) \phi_k(Y_i) - \widehat \E(\phi_k(Y)) \right \rvert \right)\\
		\leq \left( \widehat \E \left(\phi_k(Y)^2 \right) \right)^{1/2} \left( 1 + 2 Z^{1/2} \right)  e^{\alpha D(f || \hat g_t)/4} n^{-\alpha / 4}.
	\end{multline*}
	Using $(1 + 2  Z^{1/2}) e^{\alpha D(f || \hat g_t)/4}\leq Z'$ (since $Z^{1/2} \geq 1$) and summing over $k$ gives the result.
\end{proof}

The gist of CE is that $\hat \mu_{t+1}$ and $\hat \Sigma_{t+1}$ are thought of IS estimators of $\mu_{t+1} = \mu_{A_t}$ and $\Sigma_{t+1} = \Sigma_{A_t}$, which suggests to use the bound of the previous display to control them. However, a close inspection of their definitions~\eqref{eq:CE-p-mu} and~\eqref{eq:CE-Sigma} shows that $\hat \mu_{t+1}$ and $\hat \Sigma_{t+1}$ are not exactly IS estimators of $\mu_{A_t}$ and $\Sigma_{A_t}$ for two reasons:
\begin{enumerate}
	\item they are self-normalized through the estimator $\hat p_t$;
	\item they are IS estimators of $\mu_{\hat A_t}$ and $\Sigma_{\hat A_t}$, rather than  $\mu_{A_t}$ and $\Sigma_{A_t}$.
\end{enumerate}
The first point prevents from directly using the bound of the previous display. For this reason, we start by analyzing the following quantities:
\[ w_t = \frac{\hat p_t}{p_f(\hat A_t)}, \ \hat \mu'_{t+1} = \frac{1}{\npar p_f(\hat A_t)} \sum_{i=1}^\npar \ell(Y_i) \xi_{\hat A_t}(Y_i) Y_i \]
and
\[ \hat \Sigma'_{t+1} = \frac{1}{\npar p_f(\hat A_t)} \sum_{i=1}^\npar \ell(Y_i) \xi_{\hat A_t}(Y_i) (Y_i - \mu_{\hat A_t}) (Y_i - \mu_{\hat A_t})^\top \]
where here and in the sequel, $\ell = f/\hat g_t$ and, under $\widehat \P$, the $Y_i$'s are i.i.d.\ drawn according to $\hat g_t$. Then $\hat \mu'_{t+1}$ and $\hat \Sigma'_{t+1}$ are the IS estimators of $\mu_{\hat A_t}$ and $\Sigma_{\hat A_t}$, respectively, with the IS density $\hat g_t$. In particular, we can apply the previous lemma to control them, which leads to the following bounds.

\begin{lemma} \label{lemma:CD-bounds}
	With the notation introduced above, we have
	\begin{equation} \label{eq:w}
		\widehat \E \left( \left \lvert w_t - 1 \right \rvert \right) \leq \frac{Z'}{(p_f(\hat A_t))^{1/2}} n^{-\alpha/4},
	\end{equation}
	\begin{equation} \label{eq:mu'}
		\widehat \E \left( \lone{\hat \mu'_{t+1} - \mu_{\hat A_t}} \right) \leq \frac{Z'}{p_f(\hat A_t)} d n^{-\alpha/4}
	\end{equation}
	and
	\begin{equation} \label{eq:Sigma'}
		\widehat \E \left( \lone{\hat \Sigma'_{t+1} - \Sigma_{\hat A_t}} \right) \leq \frac{ (4 + 2\lVert \mu_{\hat A_t} \rVert^2) Z'}{p_f(\hat A_t)} d^2 n^{-\alpha/4}.
	\end{equation}
\end{lemma}

\begin{proof}
	Recall that $\hat p_t = \frac{1}{n} \sum_i \ell(Y_i) \xi_{\hat A_t}(Y_i)$: applying~\eqref{eq:CD+exp-tail} with $\phi = \xi_{\hat A_t}$, we obtain
	\[ \widehat \E \left( \Lone{ \hat p_t - p_f(\hat A_t)} \right) \leq Z' p_f(\hat A_t)^{1/2} n^{-\alpha/4} \]
	which gives~\eqref{eq:w} by dividing by $p_f(\hat A_t)$. For the second bound~\eqref{eq:mu'}, we use~\eqref{eq:CD+exp-tail} with $\phi(x) = \xi_{\hat A_t}(x) x$, corresponding to $\phi_k(x) = \xi_{\hat A_t}(x)x(k)$: then
	\begin{multline*}
		\widehat \E \left( \Lone{\frac{1}{\npar} \sum_{k=1}^{\npar} \ell(Y_i) \xi_{\hat A_t}(Y_i) Y_i - \widehat \E(Y \xi_{\hat A_t}(Y))} \right) \leq Z'\\\times \left ( \sum_{k=1}^d \sqrt{\widehat \E\left( Y(k)^2 \xi_{\hat A_t}(Y) \right)} \right) n^{-\alpha/4}.
	\end{multline*}
	Since $Y \sim f$, we have $\widehat \E\left( Y(k)^2 \xi_{\hat A_t}(Y) \right) \leq \widehat \E(Y(k)^2) = 1$ and so using this bound and dividing by $p_f(\hat A_t)$, we obtain
	\[ \widehat \E \left( \Lone{\frac{1}{\npar p_f(\hat A_t)} \sum_{i=1}^\npar \ell(Y_i) \xi_{\hat A_t}(Y_i) Y_i - \widehat \E \left( Y \mid Y \in \hat A_t \right)} \right) \leq \frac{Z'}{p_f(\hat A_t)} d n^{-\alpha/4}. \]
	Recalling the definitions
	\[ \hat \mu'_{t+1} = \frac{1}{\npar p_f(\hat A_t)} \sum_{i=1}^\npar \ell(Y_i) \xi_{\hat A_t}(Y_i) Y_i \ \text{ and } \ \mu_{\hat A_t} = \widehat \E \left( Y \mid Y \in \hat A_t \right) \]
	we see that this exactly~\eqref{eq:mu'}. The bound~\eqref{eq:Sigma'} for the variance follows along similar lines by considering $\phi(x) = (x - \mu_{\hat A_t}) (x - \mu_{\hat A_t})^\top \xi_{\hat A_t}(x)$. For this choice of $\phi$, starting from~\eqref{eq:CD+exp-tail} and dividing by $p_f(\hat A_t)$, we obtain similarly as above
	\begin{equation} \label{eq:bound-Sigma}
		\widehat \E \left( \Lone{ \hat \Sigma_{t+1}' - \Sigma_{\hat A_t} } \right) \leq Z' \left ( \sum_{1 \leq i, j \leq d} \sqrt{\widehat \E\left( Z_i Z_j \right)} \right) n^{-\alpha/4}
	\end{equation}
	with $Z_i = (Y(i) - \mu_{\hat A_t}(i))^2$. Since $Z_i$ and $Z_j$ are independent under $\widehat \P$ for $i \neq j$, we have
	\begin{align*}
		\sum_{i,j} \sqrt{\widehat \E(Z_i Z_j)} & = \sum_{i=1}^d \sqrt{\widehat \E(Z_i^2)} + \sum_{i\neq j} \sqrt{\widehat \E(Z_i) \widehat \E(Z_j)} \\
		& \leq \sum_{i=1}^d \sqrt{\widehat \E(Z_i^2)} + \left( \sum_{i=1}^d \sqrt{\widehat \E(Z_i)} \right)^2.
	\end{align*}
	using for the last inequality that the $Z_i$'s are non-negative. Using that $\widehat \E(Y(i)^k) = 0$ for $k = 1,3$ and $5$, that $\widehat \E(Y(k)^2) = 1$ and that $\widehat \E(Y(k)^4) = 3$ (because $Y \sim f$), we can compute (bounds on) the first and second moments of the $Z_i$'s. For the first moment, we have
	\[ \widehat \E \left( Z_i \right) = \widehat \E \left( (Y(i) - \mu_{\hat A_t}(i))^2 \right) = 1 + \mu_{\hat A_t}(i)^2 \leq 1 + \lVert \mu_{\hat A_t} \rVert^2 \]
	and for the second moment, we have
	\[ \widehat \E \left( Z^2_i \right) = \widehat \E \left( (Y(i) - \mu_{\hat A_t}(i))^4 \right) = \mu_{\hat A_t}(i)^4 + 6 \mu_{\hat A_t}(i)^2 + 3 \leq (\mu_{\hat A_t}(i)^2 + 3)^2 \]
	and so $\widehat \E \left( Z^2_i \right) \leq (\lVert \mu_{\hat A_t} \rVert^2 + 3)^2$. This gives
	\[ \sum_{i,j} \sqrt{\widehat \E(Z_i Z_j)} \leq d (\lVert \mu_{\hat A_t} \rVert^2 + 3) + d^2 (1+\lVert \mu_{\hat A_t} \rVert^2) \leq d^2 (4+2\lVert \mu_{\hat A_t} \rVert^2). \]
	Plugging in this inequality into~\eqref{eq:bound-Sigma} gives the result.
\end{proof}

\begin{corollary} \label{cor:induction-1}
	Assume that:
	\begin{itemize}
		\item the stochastic induction hypothesis holds at time $t$;
		\item $\inf_d \rho > 0$;
		\item for every $d$, $\varphi$ has no atom;
		\item $m \to \infty$;
		\item $n \gg d^{8/\alpha}$.
	\end{itemize}
	Then $\lVert \hat \mu_{t+1} \rVert$, $\Psi(\hat \Sigma_{t+1})$ and $1/p_f(\hat A_{t+1})$ are bounded whp, i.e., the stochastic induction hypothesis holds at time $t+1$.
\end{corollary}

\begin{proof}
	Stochastic induction holding at time $t$ gives $\P(\calE) \to 1$ by Lemma~\ref{lemma:alpha_*}, so we can assume without loss of generality that the event $\calE$ holds almost surely. Let us first prove that $\lVert \hat \mu_{t+1} \rVert$ is bounded whp. By the stochastic induction hypothesis, $1/p_f(\hat A_t)$ is bounded whp, so $\lVert \mu_{\hat A_t} \rVert$ and $\Psi(\Sigma_{\hat A_t})$ are bounded whp by Corollary~\ref{cor:Sigma-proj-mu-A-bounded-stoch}. Therefore, it is enough to prove that $\lVert \hat \mu_{t+1} - \mu_{\hat A_t} \rVert \Rightarrow 0$. By definition we have
	\[ \hat \mu'_{t+1} = \frac{1}{\npar p_f(\hat A_t)} \sum_{i=1}^\npar \ell(Y_i) \xi_{\hat A_t}(Y_i) Y_i = w_t \hat \mu_{t+1} \]
	and so
	\begin{align*}
		\lVert \hat \mu_{t+1} - \mu_{\hat A_t} \rVert & \leq \lVert \hat \mu_{t+1} - \hat \mu'_{t+1} \rVert + \lVert \hat \mu'_{t+1} - \mu_{\hat A_t} \rVert\\
		& = \frac{\lvert w_t - 1\rvert}{w_t} \lVert \hat \mu'_{t+1} \rVert + \lVert \hat \mu'_{t+1} - \mu_{\hat A_t} \rVert\\
		& \leq \frac{\lvert w_t - 1\rvert}{w_t} \lVert \mu_{\hat A_t} \rVert + \left(1 + \frac{\lvert w_t - 1\rvert}{w_t} \right) \lVert \hat \mu'_{t+1} - \mu_{\hat A_t} \rVert\\
		& \leq \frac{\lvert w_t - 1\rvert}{w_t} \lVert \mu_{\hat A_t} \rVert + \left(1 + \frac{\lvert w_t - 1\rvert}{w_t} \right) \lone{\hat \mu'_{t+1} - \mu_{\hat A_t}}.
	\end{align*}
	By the stochastic induction hypothesis and Lemma~\ref{lemma:E-Z}, $1/p_f(\hat A_t)$ and $Z'$ are bounded whp: therefore, we get $\widehat \E(\lvert w_t - 1 \rvert) \Rightarrow 0$ and $\widehat \E(\lone{\hat \mu'_{t+1} - \mu_{\hat A_t}}) \Rightarrow 0$ by Lemma~\ref{eq:mu'} which implies that $w_t \Rightarrow 1$ and $\lone{\hat \mu'_{t+1} - \mu_{\hat A_t}} \Rightarrow 0$. Since $\lVert \mu_{\hat A_t} \rVert$ is bounded whp, the last bound of the previous display implies that $\lVert \hat \mu_{t+1} - \mu_{\hat A_t} \rVert \Rightarrow 0$ as desired.
	
	Let us now prove that $\Psi(\hat \Sigma_{t+1})$ is bounded whp. Since $\Psi(\Sigma_{\hat A_t})$ is bounded whp,  $\lVert \Sigma_{\hat A_t} - I \rVert$, $\lambda_d(\Sigma_{\hat A_t})$ and $1/\lambda_1(\Sigma_{\hat A_t})$ are bounded whp by Lemma~\ref{lemma:equivalences}. Moreover,
	\begin{equation} \label{eq:3}
		\lVert \hat \Sigma_{t+1} - I \rVert \leq \lVert \hat \Sigma_{t+1} - \hat \Sigma'_{t+1} \rVert + \lVert \hat \Sigma'_{t+1} - \Sigma_{\hat A_t} \rVert + \lVert \Sigma_{\hat A_t} - I \rVert.
	\end{equation}
	We have just seen that the last term $\lVert \Sigma_{\hat A_t} - I \rVert$ of the right-hand side of the previous inequality is bounded whp; the second term $\lVert \hat \Sigma'_{t+1} - \Sigma_{\hat A_t} \rVert$ converges to $0$ (in distribution) because $\lVert \hat \Sigma'_{t+1} - \Sigma_{\hat A_t} \rVert \leq \lvert \hat \Sigma'_{t+1} - \Sigma_{\hat A_t} \rvert$, the latter vanishing in view of~\eqref{eq:Sigma'} (again, $Z'$, $\lVert \mu_{\hat A_t} \rVert$ and $1/p_f(\hat A_t)$ are bounded whp). Finally, the definition~\eqref{eq:CE-Sigma} of $\hat \Sigma_{t+1}$ can be rewritten as
	\[ \hat \Sigma_{t+1} = \frac{1}{\npar \hat p_t} \sum_{i=1}^\npar \ell(Y_i) \xi_{\hat A_t}(Y_i) (Y_i - \hat \mu_{t+1}) (Y_i - \hat \mu_{t+1})^\top. \]
	Recalling that $\hat p_t = \frac{1}{n} \sum_{i=1}^n \ell(Y_i) \xi_{\hat A_t}(Y_i)$ and that $\hat \mu_{t+1} = \frac{1}{n \hat p_t} \sum_{i=1}^n \ell(Y_i) \xi_{\hat A_t}(Y_i) Y_i$, we get
	\[ \frac{1}{\npar \hat p_t} \sum_{i=1}^\npar \ell(Y_i) \xi_{\hat A_t}(Y_i) (Y_i -  \mu_{\hat A_t}) = \hat \mu_{t+1} - \mu_{\hat A_t}. \]
	Starting from the previous expression of $\hat \Sigma_{t+1}$, writing $Y_i - \hat \mu_{t+1} = a + b$ with $a = Y_i - \mu_{\hat A_t}$ and $b = \mu_{\hat A_t} - \hat \mu_{t+1}$ and expanding the product, we get
    \[ \hat \Sigma_{t+1} = \frac{1}{\npar \hat p_t} \sum_{i=1}^\npar \ell(Y_i) \xi_{\hat A_t}(Y_i) (Y_i -  \mu_{\hat A_t}) (Y_i - \mu_{\hat A_t})^\top - (\mu_{\hat A_t}- \hat \mu_{t+1}) (\mu_{\hat A_t}- \hat \mu_{t+1})^\top \]
    which finally leads to
	\[ \hat \Sigma_{t+1} = \frac{1}{w_t} \hat \Sigma'_{t+1} - (\mu_{\hat A_t}- \hat \mu_{t+1}) (\mu_{\hat A_t}- \hat \mu_{t+1})^\top. \]
	Since $\lVert x x^\top \rVert = \lVert x \rVert^2$, we get
	\begin{align*}
		\lVert \hat \Sigma_{t+1} - \hat \Sigma'_{t+1} \rVert & \leq \frac{\lvert w_t - 1 \rvert}{w_t} \lVert \hat \Sigma'_{t+1} \rVert + \lVert \mu_{\hat A_t}- \hat \mu_{t+1} \rVert^2\\
		& \leq \frac{\lvert w_t - 1 \rvert}{w_t} \lVert \hat \Sigma'_{t+1} - \Sigma_{\hat A_t} \rVert + \frac{d \lvert w_t - 1 \rvert}{w_t} \lambda_d(\Sigma_{\hat A_t}) + \lVert \mu_{\hat A_t}- \hat \mu_{t+1} \rVert^2,
	\end{align*}
	using the triangle inequality and $\lVert \Sigma_{\hat A_t} \rVert \leq d \lambda_d(\Sigma_{\hat A_t})$ for the last inequality. We have argued that $\lvert w_t - 1 \rvert$, $\lVert \hat \Sigma'_{t+1} - \Sigma_{\hat A_t} \rVert$ and $\lVert \hat \mu_{t+1} - \mu_{\hat A_t} \rVert \Rightarrow 0$; moreover, $\lambda_d(\Sigma_{\hat A_t})$ and $1/w_t$ are bounded whp; finally, the convergence $\lvert w_t - 1 \rvert \Rightarrow 0$ can actually be strengthened to $d \lvert w_t - 1 \rvert \Rightarrow 0$ in view of~\eqref{eq:w}, because we have chosen $\alpha$ such that $d n^{-\alpha / 4} \to 0$. Therefore, all the terms in the upper bound of the previous display vanish, which implies that $\lVert \hat \Sigma_{t+1} - \hat \Sigma'_{t+1} \rVert \Rightarrow 0$. Going back to~\eqref{eq:3} we see that this implies that $\lVert \hat \Sigma_{t+1} - I \rVert$ is bounded whp, which directly implies that $\lambda_d(\hat \Sigma_{t+1})$ is also bounded whp since
	\[ \lVert \hat \Sigma_{t+1} - I \rVert^2 = \sum_i (\lambda_i(\hat \Sigma_{t+1})- 1)^2 \geq (\lambda_d(\hat \Sigma_{t+1})- 1)^2. \]
	Furthermore,
	\[ \lambda_1(\hat \Sigma_{t+1}) \geq \lambda_1(\Sigma_{\hat A_t}) - \lVert \hat \Sigma_{t+1} - \Sigma_{\hat A_t} \rVert \]
	by Lemma~\ref{lemma:L1-matrix}. Since $1/\lambda_1(\Sigma_{\hat A_t})$ is bounded whp and $\lVert \hat \Sigma_{t+1} - \Sigma_{\hat A_t} \rVert \Rightarrow 0$, the inequality of the previous display implies that $1/\lambda_1(\hat \Sigma_{t+1})$ is bounded whp. Thus, we have proved that $\lambda_d(\hat \Sigma_{t+1})$, $1/\lambda_1(\hat \Sigma_{t+1})$ and $\lVert \hat \Sigma_{t+1} - I \rVert$ are bounded whp, which implies that $\Psi(\hat \Sigma_{t+1})$ is bounded whp by Lemma~\ref{lemma:equivalences}. This achieves to prove that $\Psi(\hat \Sigma_{t+1})$ is bounded whp.
	
	In order to conclude the proof, it remains to prove that $1/p_f(\hat A_{t+1})$ is bounded whp. Using Corollary~\ref{cor:jensen} with $B = \hat A_{t+1}$ and $g = \hat g_{t+1}$, we obtain
	\[ p_f(\hat A_{t+1}) \geq p_{\hat g_{t+1}}(\hat A_{t+1}) \exp \left( - \Psi(\Sigma^{\hat g_{t+1}}_{\hat A_{t+1}}) - \frac{1}{2} \lVert \mu^{\hat g_{t+1}}_{\hat A_{t+1}} \rVert^2 \right). \]
	By Lemma~\ref{lemma:init-induction-hat-g_t}, $1/p_{\hat g_{t+1}}(\hat A_{t+1})$ is bounded whp, and so we only have to prove that $\Psi(\Sigma^{\hat g_{t+1}}_{\hat A_{t+1}})$ and $\lVert \mu^{\hat g_{t+1}}_{\hat A_{t+1}} \rVert$ are bounded whp. But since $\lVert \hat \mu_{t+1} \rVert$, $\hat \Psi(\Sigma_{t+1})$ and $1/p_{\hat g_{t+1}}(\hat A_{t+1})$ are bounded whp, this follows precisely from Corollary~\ref{cor:boundedness-sigma-mu} with $g = \hat g_{t+1}$ and $B = \hat A_{t+1}$.
 \end{proof}

%%%%%%%%%%%%%%%%%%%%%%%%%%%%%%%%%%%%%%%%%%%%%%
%% Single Appendix:                         %%
%%%%%%%%%%%%%%%%%%%%%%%%%%%%%%%%%%%%%%%%%%%%%%
\begin{appendix}
\section*{Proof of Proposition~\ref{pro:E-DKL}}%% if no title is needed, leave empty \section*{}.
We will first prove that 
\begin{align*}
     \mathbb{E}[D(f || \hat g_A)]
     &= D(f || g_A) \\
     &+ \frac{1}{2}\biggl[\sum_{i=1}^{d}\left(\psi\biggl(\frac{\npar-i}{2}\biggr) + \log\biggl(\frac{2}{\npar}\biggr)\right) \\
     &+ \frac{d+2}{\npar-d-2}\tr(\Sigma^{-1}_A) + \frac{d}{\npar-d-2} + \frac{d+2}{\npar-d-2}\mu_A^\top\Sigma^{-1}_A\mu_A \biggr]
\end{align*}
with $\psi$ the digamma function.
Using Lemma~\ref{lemma:D-g-mid-A-g'} with $g|_A = f$ and $g'=\hat g_A$,
\begin{align*}
    \mathbb{E}[D(f|| \hat g_A)]
    = \frac{1}{2}\biggl[\mathbb{E}(\log|\hat{\Sigma}_A|) + \mathbb{E}(\tr(\hat{\Sigma}_A^{-1})) + \mathbb{E}(\hat{\mu}_A^\top\hat{\Sigma}_A^{-1}\hat{\mu}_A) - d \biggr].
\end{align*}
According to \cite[pg40 and 108]{chatfield_introduction_1980}, the law of $\npar\hat{\Sigma}_A$ is the Wishart distribution with the parameters $\Sigma_A$ et $(\npar-1)$: $W_d(\Sigma_A,\npar-1)$. From \cite{bishop_pattern_2006}, 
\begin{align*}
    &\mathbb{E}(\npar\hat{\Sigma}_A) = (\npar-1)\Sigma_A \text{ and }\\
    &\mathbb{E}(\log|\npar\hat{\Sigma}_A|) = \sum_{i=1}^d\biggl(\psi\biggl(\frac{\npar-1}{2} +\frac{1-i}{2}\biggr) \biggr)+ d\log(2) + \log|\Sigma_A|). 
\end{align*}
Moreover, the law of $\frac{1}{\npar}\hat{\Sigma}_A^{-1}$ is the inverse-Wishart distribution with  parameters $\Sigma_A^{-1}$ et $(\npar-1)$: $W_d^{-1}(\Sigma_A^{-1},\npar-1)$ \cite{mardia_multivariate_2006}. We have
\begin{align*}
    &\mathbb{E}\biggl(\frac{1}{\npar}\hat{\Sigma}_A^{-1}\biggr) = \frac{1}{(\npar-1)-d-1}\Sigma_A^{-1}, \\
    &\mathbb{E}(\log|\hat{\Sigma}_A|) = \mathbb{E}(\log|\npar\hat{\Sigma}_A|) - \mathbb{E}(\log(\npar)^d)
        = \sum_{i=1}^d\biggl(\psi\biggl(\frac{\npar-i}{2}\biggr) + d\log\biggl(\frac{2}{\npar}\biggr)\biggr) + \log|\Sigma_A| \\
    &\text{and } \mathbb{E}(\text{tr}(\hat{\Sigma}_A^{-1})) 
    = \text{tr}(\mathbb{E}(\hat{\Sigma}_A^{-1})) 
    = \frac{\npar}{\npar-d-2}\text{tr}(\Sigma_A^{-1}).
\end{align*}
    Since $\hat{\mu}_A$ and $\hat\Sigma_A$ are the sample mean and sample covariance matrix of normally distributed samples respectively, they are independent, so 
    \begin{align*}
        \mathbb{E}(\hat{\mu}_A^\top\hat{\Sigma}_A^{-1}\hat{\mu}_A) = \text{tr}(\mathbb{E}(\hat{\Sigma}_A^{-1}\hat{\mu}_A\hat{\mu}_A^\top))
        = \text{tr}(\mathbb{E}(\hat{\Sigma}_A^{-1})\,\mathbb{E}(\hat{\mu}_A\hat{\mu}_A^\top)) 
        \\= \text{tr}\left(\frac{\npar}{\npar-d-2}\Sigma_A^{-1}\mathbb{E}(\hat{\mu}_A\hat{\mu}_A^\top)\right).
    \end{align*}
    From the equality $\mathbb{E}((\hat\mu_A-\mu_A)(\hat\mu_A -\mu_A)^\top) = \mathbb{E}(\hat{\mu}_A\hat{\mu}_A^\top)-\mu_A\mu_A^\top$, and by denoting \\ $\hat S = \frac{1}{\npar}\sum_{k=1}^{\npar}(Y_{A,k} - \mu_A)(Y_{A,k} - \mu_A)^\top$, it can be shown that 
    \begin{align*}
        \mathbb{E}(\hat{\mu}_A\hat{\mu}_A^\top) = \frac{1}{\npar}\mathbb{E}(\hat S) + \mu_A \mu_A^\top
    \end{align*}
    Since $nS \sim W_d(\Sigma,n)$, we have 
    \[\mathbb{E}(\hat{\mu}_A^\top\hat{\Sigma}_A^{-1}\hat{\mu}_A) = \frac{\npar}{\npar-d-2}\frac{d}{\npar} + \frac{\npar}{\npar-d-2}\mu_A^\top\Sigma_A^{-1}\mu_A \]
    Assembling the previous expressions gives the announced expression of $\mathbb{E}[D(f || \hat g_A)]$. Let us now discuss the how each term scales with $d$.
    The digamma function $\psi$ has the following bounds \cite{alzer_inequalities_1997}:
    \begin{align*}
    \forall x > 0, \, \log x - \frac{1}{x} \leq \psi(x) \leq \log x - \frac{1}{2x}    
    \end{align*}
    We have then 
    \begin{align*}
        K + K' \leq \sum_{i=1}^d\biggl(\log\biggl(\frac{\npar}{2}\biggr)-\psi\left(\frac{\npar-i}{2}\right)\biggr) \leq K + 2K'
    \end{align*}
    with 
    \begin{align*}
        K = \sum_{i=1}^d \biggl(\log\biggl(\frac{\npar}{2}\biggr)-\log\biggl(\frac{\npar-i}{2} \biggr)\biggr) = \sum_{i=1}^d\biggl(-\log\biggl(1 - \frac{i}{\npar}\biggr)\biggr)
    \end{align*}
    and
    \begin{align*}
        K' = \sum_{i=1}^d\frac{1}{\npar-i}
    \end{align*}
    In the case $\npar \gg d$, $K = \frac{d^2}{2\npar} + o\left(\frac{d^2}{\npar}\right)$ and $K' = \frac{d}{\npar} + o\left(\frac{d}{\npar}\right) = o\left(\frac{d^2}{\npar}\right)$. So, \[\sum_{i=1}^{d}\left(\psi\biggl(\frac{\npar-i}{2}\biggr) + \log\biggl(\frac{2}{\npar}\biggr)\right) = -\frac{d^2}{2\npar} + o\left(\frac{d^2}{\npar}\right).\]
    Moreover, since \[ \frac{d}{\lambda_d(\Sigma_A)} \leq \tr(\Sigma_A^{-1}) \leq \frac{d}{\lambda_1(\Sigma_A)} \] and that $1/\lambda_1(\Sigma_A)$ is bounded by Corollary~\ref{cor:Sigma-proj-mu-A-bounded}, there exists $C > 0$ such that $\tr(\Sigma_A^{-1}) = Cd + o(d)$. In addition, $\mu_A^\top\Sigma_A^{-1}\mu_A \leq \lVert \mu_A \rVert^2/\lambda_1(\Sigma_A)$ which is bounded by the same Corollary. Therefore,
    \begin{align*}
        \frac{d+2}{\npar-d-2}\text{tr}(\Sigma^{-1}) &= C\frac{d^2}{\npar} + o\left(\frac{d^2}{\npar}\right), \\ 
        \frac{d}{\npar-d-2} &= \frac{d}{\npar} + o\left(\frac{d}{\npar}\right) = o\left(\frac{d^2}{\npar}\right),\\ \text{and }\frac{d+2}{\npar-d-2}\mu_A^\top\Sigma_A^{-1}\mu_A  &= \mu_A^\top\Sigma_A^{-1}\mu_A \left(\frac{d}{\npar}\right) + o\left(\frac{d}{\npar}\right) = o\left(\frac{d^2}{\npar}\right).
    \end{align*}
    Therefore,
    \[
     \mathbb{E}[D(f || \hat g_A)]
     = D(f || g_A) 
     + \frac{1}{2}\left(C-\frac{1}{2}\right)\frac{d^2}{\npar} + o\left(\frac{d^2}{\npar}\right).
    \]
    So $\sup_d \E(D(f || \hat g_A)) < \infty$ if $\npar \gg d^2$, and $\E(D(f || \hat g_A)) \to \infty$ if $\npar \ll d^2$.
\end{appendix}
%%%%%%%%%%%%%%%%%%%%%%%%%%%%%%%%%%%%%%%%%%%%%%
%% Multiple Appendixes:                     %%
%%%%%%%%%%%%%%%%%%%%%%%%%%%%%%%%%%%%%%%%%%%%%%
%\begin{appendix}
%\section{???}
%
%\section{???}
%
%\end{appendix}

%%%%%%%%%%%%%%%%%%%%%%%%%%%%%%%%%%%%%%%%%%%%%%
%% Support information, if any,             %%
%% should be provided in the                %%
%% Acknowledgements section.                %%
%%%%%%%%%%%%%%%%%%%%%%%%%%%%%%%%%%%%%%%%%%%%%%
\begin{acks}[Acknowledgements]
The first author J.\ Beh is enrolled in a Ph.D. program co-funded by ONERA - The French Aerospace Lab and the University Research School EUR-MINT (State support managed by the National Research Agency for Future Investments program bearing the reference ANR-18-EURE-0023). Their financial supports are gratefully acknowledged. The authors would also like to thank J\'er\^ome Morio for his precious support and his feedback on a preliminary version of the paper, as well as the reviewers whose feedback helped improve the overall quality of the paper.
\end{acks}

\end{document}